\documentclass[11pt]{amsart} 
\usepackage{amsfonts,newlfont,latexsym}
\usepackage{euscript}
\usepackage{amssymb}
\usepackage{amsmath}
\usepackage{color}
\usepackage{hyperref}
\definecolor{orange}{rgb}{1,0.5,0}

\numberwithin{equation}{section}

      \newtheorem{Theorem}{Theorem}[section]
      
      \newtheorem{lemma}[Theorem]{Lemma}
      \newtheorem{Lemma}[Theorem]{Lemma}
      \newtheorem{corollary}[Theorem]{Corollary}
      
      \newtheorem{prop}[Theorem]{Proposition}
      \newtheorem{remark}[Theorem]{Remark}

\newtheorem{Proposition}[Theorem]{Proposition}

\newcommand{\R}{\mathbb R}

\newcommand{\Z}{\mathbb Z}

\def \e{\varepsilon}

\def\Proof{{\em Proof}\,: }
\def\proof{{\em Proof}\,: }
\def\QED{~\hfill~ $\diamond$ \vspace{7mm}}

\def \R{{\mathbb R}}
\def \Z{{\mathbb Z}}
\def \N{{\mathbb N}}

\def\holder{H\"{o}lder }
\def \smooth{$C^{\infty}$ }

\definecolor{rjs}{rgb}{.1,.4,.7}

\author[Alexander Gorodnik and Ralf Spatzier]
{Alexander Gorodnik$^\ast$ and  Ralf Spatzier$^{\ast \ast}$}

\title[Exponential Mixing of Nilmanifold Automorphisms]
{Exponential Mixing of Nilmanifold Automorphisms}

\thanks{$^\ast$ Supported in part by EPSRC grant EP/H000091/1 and ERC grant 239606}
\thanks{$^{\ast \ast}$ Supported in part by NSF grant DMS-0906085}

 \address{Department of Mathematics, University of Bristol, Bristol, BS8 1TW, U.K. }

\email{a.gorodnik@bristol.ac.uk}

 \address{Department of Mathematics, University of Michigan, Ann Arbor, MI 48109.}

\email{spatzier@umich.edu}

\begin{document}

\begin{abstract}
We study dynamical properties of automorphisms of compact nilmanifolds and prove
that every ergodic automorphism is exponentially mixing and exponentially mixing of higher orders.
This allows to establish probabilistic limit theorems and regularity of solutions of 
the cohomological equation for such automorphisms. Our method is based on the quantitative equidistribution results for
polynomial maps
% established by Green and Tao 
combined with Diophantine estimates.
\end{abstract}

\maketitle

\section{Introduction}

 Dynamics and ergodic theory of toral automorphisms have been well understood   for quite some  time.
 Ergodic toral automorphisms are always  mixing and even Bernoulli \cite{Katz-bernoulli}, and 
 have dense sets of periodic points \cite{Marcus}. However, unless they are hyperbolic, the toral
 automorphisms lack the specification
 property and, in particular, don't have Markov partitions \cite{l}. Nonetheless, it is known that
 ergodic toral automorphisms satisfy the central limit theorem and its refinements \cite{Leonov,LB99}.
 Regarding the quantitative aspects Lind established exponential mixing for ergodic toral
 automorphisms using Fourier analysis \cite{Lind}. 
Surprisingly, some of these ergodic-theoretic properties turned out to be more delicate
for automorphisms of compact nilmanifolds and still remained unexplored.
 In particular, the exponential mixing, which is one of the main
 results of this paper,  has not been established and does not easily follow using the harmonic analysis
 on nilpotent Lie groups.
% which is the natural substitute of the Fourier analysis in this set-up.
 %Instead, we  will use the work of Green and Tao on equidistribution properties of polynomial sequences in
 %nilmanifolds \cite{Green-Tao}.  They essentially provide  a dichotomy: either the sequence is
 %equidistributed or the sequence is well approximated by rationals in a suitable sense.     

\subsection{Exponential mixing}
Let $G$ be a simply connected nilpotent Lie group and $\Lambda$ a discrete 
cocompact subgroup.  The space $X=G/\Lambda$ is called a {\em  compact nilmanifold}.  
%More generally, one considers {\em \compact  infra-nilmanifolds} $X$ which are quotients of $N$ by a
%lattice in $N \rtimes K$ where $K$ is a compact subgroup of the automorphism group of $N$.
%Infra-nilmanifolds are always finitely covered by nilmanifolds \cite{Lee}. 
An {\em automorphism} $\alpha$ of $X$ is a diffeomorphism of $X$ which lifts to an automorphism of $G$.
We denote by $\mu$ the Haar probability measure on $X$. Then 
$\alpha$ preserves $\mu$. The ergodic-theoretic properties of the dynamical system $\alpha\curvearrowright
(X,\mu)$ have been studied by Parry \cite{Parry}.
He proved that an automorphism is ergodic
  if and only if the induced map on the maximal toral quotient is ergodic,
and every ergodic automorphism satisfies the Kolmogorov property.
In particular, it is mixing of all orders. In this paper we establish quantitative mixing properties
of such automorphisms. We fix a right-invariant Riemannian metric on $G$ which also defines 
a metric on $X$ and denote by $C^\theta(X)$ the space of $\theta$-\holder functions on $X$. 

Now we state the first main result of the paper.

\begin{Theorem}\label{exp mixing}
Let $\alpha$ be an ergodic automorphism of a compact nilmanifold $X$.
Then there exists $\rho=\rho(\theta)\in (0,1)$ such that for all 
$f_0,f_1\in C^\theta(X)$ and $n\in \mathbb{N}$,
\begin{align*}
 \int_X f_0(x)f_1(\alpha^n(x))\,d\mu(x)=\left(\int _X {f}_0\,  d \mu\right)\left(\int _X f_1\, d
   \mu\right)
+O\left(\rho ^n \|f_0 \|_{C^\theta}\|f_1\|_{C^\theta}  \right).
\end{align*}
\end{Theorem} 

The proof of Theorem \ref{exp mixing} is based on an equidistribution result  for the exponential map
established in Section \ref{box maps} (see Corollary
\ref{cor:line} below), which is deduced from the work of Green and Tao \cite{Green-Tao}.
This result shows that images of boxed under the exponential map are equidistributed in $X$ provided that
a certain Diophantine condition holds.
We complete the proof of Theorem \ref{exp mixing} in Section \ref{sec:mix}.  The main idea is to relate
the correlations $\left<f_0,f_1\circ\alpha^n\right>$ to averages along suitable foliations in $X$
and apply the equidistribution result established in Section \ref{box maps}.
In order to verify the Diophantine condition we use the Diophantine properties of algebraic numbers.
This leads to the proof of Theorem \ref{exp mixing} under an irreducibility condition,
and the proof of the theorem in general uses an inductive argument.
%Let us note that the idea using averages over foliations to show mixing properties had been used
%before, e.g. in \cite{Rudolph}.  Dolgopyat uses this idea in his work on central limit theorems for
%partially hyperbolic systems \cite{Dolgopyat}.

%We get other applications of this approach.  First, in Section \ref{sec-unstable}, we prove equidistribution results for  unstable manifolds, Theorem \ref{th:unstable}.   We  later use this result  in our applications to central limit theorems for ergodic automorphisms of nilmanifolds in Section \ref{CLT}. 

We also establish multiple exponential mixing for ergodic automorphisms of compact nilmanifolds.
% Let us note that multiple mixing results are never amenable to an approach via representation theory.  
For ergodic toral automorphisms, multiple exponential mixing was proved by  P\`{e}ne \cite{Pene} and Dolgopyat \cite{Dolgopyat}.

\begin{Theorem}\label{thm_mixing_lipschitz}
Let $\alpha$ be an ergodic automorphism of a compact nilmanifolds $X$.
Then there exists $\rho=\rho(\theta)\in (0,1)$ such that for all 
$f_0,\ldots,f_s\in C^\theta(X)$ and $n_0,\ldots,n_s\in \mathbb{N}$,
\begin{align*}
 \int_X \left(\prod_{i=0}^s f_i(\alpha^{n_i}(x))\right)\,d\mu(x)=
\prod_{i=0}^s \left(\int _X {f}_i\,  d \mu\right)
+O\left(\rho ^{\min_{i\ne j} |n_i-n_j|} \prod_{i=0}^s\|f_i \|_{C^\theta}  \right).
\end{align*}
\end{Theorem} 

The proof of Theorem \ref{thm_mixing_lipschitz} is given in Section \ref{sec:multiple}.
The first step of the proof is to establish an equidistribution result for 
images of exponential map in $X\times\cdots \times X$
(see Proposition \ref{p:mmult}). Then we approximate higher order correlations
by averages of the exponential map. As in the proof of Theorem \ref{exp mixing}, 
we first consider the irreducible case and then deduce the theorem in general
using an inductive argument.

\subsection{Probabilistic limit theorems} 
It is well-known that the exponential mixing property is closely related to other chaotic properties
of dynamical systems and, in particular, to the central limit theorem for observables $f\circ\alpha^n$.
While one does not imply the other in general,
the martingale differences approach \cite[Ch.~5]{HH} usually allows to deduce
the proof of the central limit theorem from quantitative equidistribution of unstable foliations.
Using this approach, the central limit theorem and its generalisations have been
established for ergodic toral automorphisms in \cite{Leonov, LB99}
and for ergodic automorphisms of 3-dimensional nilmanifolds in \cite{CB}.
Here we extend these results to general nilmanifolds.

\begin{Theorem}\label{th:clt}
Let $\alpha$ be an ergodic automorphism of a compact nilmanifolds $X$
and $f\in C^\theta(X)$ with $\int_X f\, d\mu=0$
which is not a measurable coboundary (i.e., $f\ne \phi\circ \alpha-\phi$
for any measurable function $\phi$ on $X$). Then there exists $\sigma=\sigma(f)>0$ such that
$$
\mu\left(\left\{x\in X:\, \frac{1}{\sqrt{n}} \sum_{i=0}^{n-1} f(\alpha^i(x))\in (a,b)
  \right\}\right) \to \frac{1}{\sqrt{2\pi}\sigma}\int_a^b e^{-{t^2}/(2\sigma^2)}\, dt
$$
as $n\to \infty$.
\end{Theorem}

We also prove 
the {\it central limit theorem for subsequences}, and the 
{\it Donsker and Strassen invariance principles} for ergodic automorphisms of nilmanifolds.  
We refer to Section  \ref{CLT} for a detailed discussion of the results. 
The main ingredient of the proof is the exponential equidistribution of leaves of
unstable foliations, which is established for this purpose in Section \ref{sec-unstable}.

\subsection{Cohomological equation}     
Let $\alpha$ be a measure-preserving transformation of a probability space $(X, \mu)$
and $f:X\to \mathbb{R}$ is a measurable function. The functional equation
\begin{equation}\label{eq:cobb}
f=\phi\circ \alpha -\phi, \quad\quad \phi:X\to \mathbb{R},
\end{equation}
is called the {\it cohomological equation}. This equation plays important role in many
aspects of the theory of dynamical systems (for instance,
existence of smooth invariant measures, existence of conjugacies, 
existence of isospectral deformations, rigidity of group actions). 
If a measurable solution $\phi$ of (\ref{eq:cobb}) exists, the function $f$
is called a {\it measurable coboundary}.
It is easy to see that a solution of (\ref{eq:cobb}) is unique 
(up to measure zero) up to an additive constant when $\alpha$
is ergodic with respect to $\mu$. 

We will apply the exponential mixing property to investigate regularity of solutions of
the cohomological equation.
  
\begin{Theorem}   \label{th:cobb}
Let $\alpha$ be an ergodic automorphism  of a compact nilmanifold $X$,
and let $f\in C^{\infty}(X)$ be such that (\ref{eq:cobb}) has a measurable solution.
Then there exists a $C^\infty$ solution of (\ref{eq:cobb}). 
\end{Theorem}

The first result of this type was proved by Livsic  \cite{Livsic} for Anosov diffeomorphism and flows. 
% in what is called the measurable Livsic' theorem.  
More precisely, if $\alpha$ is an Anosov diffeomorphism and the given $C^\infty$ function $f$
is a measurable coboundary, then the cohomological equation (\ref{eq:cobb}) has a $C^\infty$ solution $\phi$. 
There are also  versions of this result for H\"older functions and $C^k$ functions.
Recently, Wilkinson \cite{Wilkinson} has generalised Livsic' results to
partially hyperbolic diffeomorphisms that satisfy the accessibility property.
Automorphisms of nilmanifolds however do not have the accessibility property.
In fact, the problem of regularity of solutions of the coboundary equation for ergodic toral
automorphisms, which are not hyperbolic,
turns out to be quite subtle \cite{Veech,LB98}. 
Veech \cite{Veech} has constructed an example of $f\in C^1(\mathbb{T}^d)$ 
which sums to zero along periodic orbits, but the cohomological equation (\ref{eq:cobb})
has no $C^1$ solutions.
By \cite{Veech}, if $f\in C^k(\mathbb{T}^d)$ with $k>d$ and (\ref{eq:cobb})
has a measurable solution, then there exists a solution in $C^{k-d}(\mathbb{T}^d)$. 
We are not aware of any results regarding regularity of solutions of (\ref{eq:cobb})
for a general  ergodic toral automorphism when $f\in C^k(\mathbb{T}^d)$ with $k<d$.

Theorem \ref{th:cobb} is proved in Section \ref{sec:cohomology}.
We use a construction from Section \ref{CLT} to show that
there exists a square-integrable solution.
Then we use a new method of proving smoothness as developed by Fisher, Kalinin and Spatzier in \cite{FKS}:
we consider the solution as a distribution on the space of \holder functions and
study its regularity along the stable, unstable and central foliations
of $\alpha$. While regularity along the first two foliations can be deduced using
the standard contraction argument, the regularity along the central foliation 
is deduced from the exponential mixing property.

\subsection{Further generalisations}
\begin{itemize}
\item
We note that the results established here can be generalised to {\it affine diffeomorphisms} of a compact nilmanifold $X=G/\Lambda$.
Those are diffeomorphisms $\sigma:X \to X$ that can be lifted to affine maps $\tilde\sigma$ of $G$,
i.e., maps $\tilde{\sigma}$ that have constant derivatives with respect to a right invariant framing of $G$. 
Since every such diffeomorphism $\sigma$ is of the form $\sigma(x)=g_0\alpha(x)$ for $g_0\in G$ and an automorphism
$\alpha$ of $X$, our method applies to such maps as well (see Remark \ref{rem:affine} below).

\item More generally, one may consider {\it infra-nilmanifolds} \cite{dekimpe0}.
Let $G$ be a simply connected nilpotent Lie group, $C$ a compact 
subgroup of $\hbox{Aut}(G)$, and $\Gamma$ a discrete torsion-free subgroup
of $G\rtimes C$ such that $G/\Gamma$ is compact. The space $Y=G/\Gamma$
is called an infra-nilmanifold. 
By \cite[Th.~1]{aus}, the group $\Lambda=G\cap \Gamma$ has finite index in $\Gamma$.
Hence, the infra-nilmanifold $Y$ is finitely covered by the nilmanifold $X=G/\Lambda$.
An affine diffeomorphism of $Y$ is a diffeomorphism which lifts
to an affine map of $G$. Every such diffeomorphism is of the form $g\mapsto g_0\alpha(g)$, where
$g_0\in G$ and $\alpha$ is an automorphism of $G$ that preserves the orbits of $\Gamma$.
By \cite[Theorem~3.4]{dekimpe}, we must have $\alpha\Gamma\alpha^{-1}=\Gamma$.
Since by \cite[Prop.~2]{aus} $\Lambda$ is the maximal normal nilpotent subgroup of $\Gamma$,
we deduce that $\alpha(\Lambda)=\alpha\Lambda\alpha^{-1}=\Lambda$. 
Therefore, every affine diffeomorphism of $Y$ lifts to an affine diffeomorphism of $X$,
and our results can be generalised to this setting.

\item Our techniques also allow to establish exponential mixing properties for {\it $\Z^k$-actions} by automorphisms
of nilmanifolds when $k\ge 2$. Since this requires more delicate Diophantine estimates, we pursue this
in a sequel paper \cite{GorodnikSpatzierII}.  This result
has found a striking application to the problem of global rigidity of smooth actions.
Given any $C^\infty$-action of $\mathbb{Z}^k$, $k\ge 2$, on a nilmanifold
that has sufficiently many Anosov elements, 
 Fisher, Kalinin and the second author showed in \cite{FKS} that this action is 
$C^{\infty}$-conjugate to an affine action on the nilmanifold.

\item
In view of the works of Katznelson \cite{Katz-bernoulli} and Parry \cite{Parry},
it is natural to ask whether ergodic automorphisms of compact nilmanifolds are {\it Bernoulli}.
Surprisingly, we could not find this result in the literature, and in Section \ref{s:ber}
we establish the Bernoulli property.
While this  easily  follows from the works of Marcuard \cite{Marcuard} and Rudolph \cite{Rudolph78}, 
and the proof does not rely on the main ideas of this paper, we
 include  this result in Section \ref{s:ber} to 
complete our discussion of ergodic properties of nilmanifold automorphisms.

\end{itemize}

\subsection*{Acknowledgements}
We are indebted to J. Rauch for discussions concerning his regularity theorems with M. Taylor.
Also we thank F. Ledrappier  for the reference to Le Borgne's work which was crucial to our applications
to the central limit theorem. A.G. would like to thank the University of Michigan for hospitality
during his visit when the work on this project had started.
R.S. thanks the University of Bristol for hospitality and support during this work.

%%%%%%%%%%%%%%%%%%%%%%%%%%%%%%%%%%
%%%%%%%%%%%%   EQUIDISTRIBUTION OF BOX MAPS  %%%%%%%%%%%

\section{Equidistribution of box maps}  \label{box maps}

Let $G$ be a simply connected nilpotent Lie group, $\Lambda$ a discrete cocompact subgroup,
and $X=G/\Lambda$ the corresponding nilmanifold equipped with the Haar probability measure $\mu$. 
We fix a a right invariant Riemannian metric $d$ on $G$ which also defines a metric on $X$.
Let $\mathcal{L}(G)$ be the Lie algebra of $G$ and $\exp:\mathcal{L}(G)\to G$ the exponential map. 
The aim of this section is to investigate distribution of images of the maps
$$
\mathbb{R}^k\to X:t\mapsto g_1\exp(\iota(t))g_2\Lambda
$$
with $g_1,g_2\in G$ and an affine map $\iota:\mathbb{R}^k\to\mathcal{L}(G)$.

The lattice subgroup $\Lambda$ defines a rational structure on $\mathcal{L}(G)$.
% we denote by $\mathcal{L}(G)_{\mathbb{Q}}\subset \mathcal{L}(G)$ the Lie algebra over $\mathbb{Q}$. 
Let $\pi:G\to G/G'$ denote the factor map, where $G'$ is the commutator subgroup.
We also have the corresponding map $D\pi:\mathcal{L}(G)\to \mathcal{L}(G/G')$.
We fix an identification $G/G'\simeq \mathcal{L}(G/G')\simeq \mathbb{R}^l$ that respects the rational structures.

We call a {\it box map} an affine map 
$$
\iota:B:=[0,T_1]\times\cdots\times [0,T_k]\to \mathcal{L}(G)
$$ of the form
\begin{equation}\label{eq:box}
\iota: (t_1,\ldots,t_k)\mapsto v+t_1w_1+\cdots +t_kw_k
\end{equation}
with $v,w_1,\ldots,w_k\in \mathcal{L}(G)$.
We denote by 
$$
|B|:=T_1\cdots T_k
$$
the volume of the box 
$B$ and by 
$$
\min(B):=\min_{i=1,\ldots,k} T_i,
$$
the length of the shortest side of $B$.

\begin{Theorem}\label{th:green-tao}
There exist $L_1,L_2>0$ such that 
for every $\delta\in (0,1/2)$ and
every box map $\iota:B\to\mathcal{L}(G)$ as in  (\ref{eq:box}), % with $T_1,\ldots, T_k\ge 1$,
 one of the following holds:
\begin{enumerate}
\item[(i)] For every Lipschitz function $f:X\to \mathbb{R}$, $u\in \mathcal{L}(G)$, and $g\in G$, 
\begin{equation}\label{eq:pos1}
\left|\frac{1}{|B|}\int_B f(\exp(u)\exp(\iota(t))g\Lambda)\, dt-\int_X f\, d\mu\right|\le \delta \|f\|_{ Lip}.
\end{equation}
\item[(ii)] There exists $z\in \mathbb{Z}^l\backslash \{0\}$ such that
\begin{equation}\label{eq:pos2}
\|z\|\ll \delta^{-L_1}\quad\hbox{and}\quad |\left<z,D\pi(w_i)\right>|\ll \delta^{-L_2}/T_i \quad \hbox{for
  all $i=1,\ldots,k$.}
\end{equation}
\end{enumerate}
\end{Theorem}

Here and in the rest of the paper we explicitly list dependences of implied constants on 
relevant parameters. In particular, in (\ref{eq:pos2}) the implied constants are independent
of the box map. 

\proof
We suppose that (i) fails for some Lipschitz function $f$,
$u\in \mathcal{L}(G)$, and $g\in G$.
Then will show that (ii) holds. We pick $L\ge 2$ such that
\begin{equation}\label{eq:initial}
\max\{\|u\|,\|v\|,T_1\|w_1\|,\ldots,T_k\|w_k\|\}\le \delta^{-L}.
\end{equation}
Making a linear change of variables in the integral~(\ref{eq:pos1}), we arrange that  $T_i\ge 1$ and $\|w_i\|\le 1$.

For $x_1,x_2,x_3\in \mathcal{L}(G)$, we consider the map 
$$
P(x_1,x_2,x_3):=\exp(x_1)\exp(x_2+x_3)\exp(-x_2)\exp(-x_1).
$$
We note that $G$ can be equipped with a structure of algebraic group so that $\exp$ 
is a polynomial isomorphism. Hence, the map $P$ can be written as
$$
P(x_1,x_2,x_3)= \exp(p_1(x_1,x_2,x_3)e_1+\cdots+ p_d(x_1,x_2,x_3)e_d)
$$
for some polynomials $p_i$.
Since $P(x_1,x_2,0)=e$, these polynomials satisfy $p_i(x_1,x_2,0)=0$.
Hence, assuming that $\|x_3\|\le 1$, we obtain
$$
|p_i(x_1,x_2,x_3)|\ll (1+\|x_1\|)^{\deg(p_i)} (1+\|x_2\|)^{\deg(p_i)}\|x_3\|,\quad i=1,\ldots,d.
$$ 
Since in the neighborhood of the origin, 
$$
d(e,P(x_1,x_2,x_3))\ll \max_{i=1,\ldots, d} |p_i(x_1,x_2,x_3)|,
$$
we deduce that there exists $C_0\ge 2$ such that
for every $\epsilon\in (0,1/2)$ and  $x_1,x_2,x_3\in\mathcal{L}(G)$ satisfying
$\|x_1\|,\|x_2\|\le (k+1)\epsilon^{-1}$ and $\|x_3\|\le   k\epsilon^{C_0}$, we have 
\begin{equation}\label{eq:exp}
d(e,P(x_1,x_2,x_3))\le \epsilon.
\end{equation}
We set $s=\lceil \delta^{-C L}\rceil$, where $C\ge C_0$ 
is sufficiently large and will be specified later
(see (\ref{eq:c_0}) and (\ref{eq:e1})--(\ref{eq:e2}) below).
Let
$$
\mathcal{N}:=\{(n_1,\ldots,n_k):\, n_i=0,\ldots, N_i-1\},
$$
where $N_i:=\lceil T_i s\rceil\ge s$. 
We consider the polynomial map
$$
p(n):=\exp(u)\exp\left(v+\sum_{i=1}^k\frac{n_i}{s}w_i\right)g,\quad n\in\mathcal{N}.
$$
For $t_i\in [\frac{n_i}{s},\frac{n_i+1}{s}]$, we apply (\ref{eq:exp}) with 
$$
x_1:=u,\;\;\; x_2:=v+\sum_{i=1}^k\frac{n_i}{s}w_i,\;\;\;
x_3:=\sum_{i=1}^k\left(t_i-\frac{n_i}{s}\right)w_i,\;\;\; \epsilon=\delta^{L}.
$$
It follows from (\ref{eq:initial}) that
\begin{align*}
\|x_1\| &\le \delta^L,\\
\|x_2\|&\le \delta^L +\sum_{i=1}^k (N_i-1)s^{-1}\|w_i\|\le \delta^L +\sum_{i=1}^k T_i\|w_i\|\le (k+1)\delta^L,\\
\|x_3\|&\le \sum_{i=1}^k s^{-1}\|w_i\|\le k s^{-1}\le  k\delta^{CL}.
\end{align*}
Hence, (\ref{eq:exp}) gives
\begin{align}\label{eq:est_m}
&d\left(p(n)\Lambda, \exp(u)\exp\left(v+\sum_{i=1}^kt_iw_i\right)g\Lambda\right)\\
\le& d\left(e, \exp(u)\exp\left(v+\sum_{i=1}^kt_iw_i\right)g p(n)^{-1}\right) \nonumber\\
= & d\left(e, \exp(u)\exp\left(v+\sum_{i=1}^kt_iw_i\right)
\exp\left(v+\sum_{i=1}^k\frac{n_i}{s}w_i\right)^{-1}\exp(u)^{-1}\right)
\le \delta^{L}. \nonumber
\end{align}
For $n=(n_1,\cdots,n_k)\in\mathcal{N}$, we set
$$
B_n:=\left[\frac{n_1}{s},\frac{n_1+1}{s}\right]\times\cdots\times \left[\frac{n_k}{s},\frac{n_k+1}{s}\right].
$$
It follows from (\ref{eq:est_m}) that for every Lipschitz function $f$ and $n\in\mathcal{N}$,
\begin{align*}
\left| f(p(n)\Lambda)|B_n|-\int_{B_n}f(\exp(u)\exp(\iota(t))g\Lambda)\,dt\right|
\le \delta^{L}s^{-k}\|f\|_{Lip}.
\end{align*}
We also observe that $B\supset \cup_{n\in\mathcal{N}} B_n$, and
$$
\left|B-\left(\bigcup_{n\in\mathcal{N}} B_n\right)\right|\le k s^{-1}T_1\cdots T_k\le ks^{-k-1} N_1\cdots N_k.
$$
Therefore, we deduce that
\begin{align*}
&\left| \sum_{n\in\mathcal{N}} f(p(n)\Lambda)|B_n|-\int_B f(\exp(u)\exp(\iota(t))g\Lambda)\,dt\right|\\
\le& \sum_{n\in\mathcal{N}} \left| f(p(n)\Lambda)|B_n|- \int_{B_n}
  f(\exp(u)\exp(\iota(t))g\Lambda)\,dt\right| +k s^{-k-1} N_1\cdots N_k\|f\|_{Lip}\\
\le & \left(\delta^{L}+k s^{-1}\right) s^{-k}N_1\cdots N_k\|f\|_{Lip},
\end{align*}
and
\begin{align*}
&\left| \frac{1}{N_1\cdots N_k}\sum_{n\in\mathcal{N}} f(p(n)\Lambda)-\frac{1}{|B|}\int_{B}
f(\exp(u)\exp(\iota(t))g\Lambda)\,dt\right|\\
\le& \left| \frac{1}{N_1\cdots N_k}\sum_{n\in\mathcal{N}} f(p(n)\Lambda)-\frac{s^k}{N_1\cdots N_k}\int_B
  f(\exp(u)\exp(\iota(t))g\Lambda)\,dt\right|\\
&+\left(\frac{1}{|B|}-\frac{s^k}{N_1\cdots N_k}\right)|B| \|f\|_{Lip}\\
\le & \left(\delta^{L}+k s^{-1}\right)\|f\|_{Lip}
+\left(1-\frac{s^k T_1\cdots T_k}{N_1\cdots N_k}\right) \|f\|_{Lip}\\
\le & \left(\delta^{L}+k s^{-1}\right)\|f\|_{Lip}
+\left(1-\frac{(N_1-1)\cdots (N_k-1)}{N_1\cdots N_k}\right) \|f\|_{Lip}\\
\le & \left(\delta^{L}+c_k s^{-1} \right)\|f\|_{Lip}\le  (\delta^{L}+c_k\delta^{CL})\|f\|_{Lip}
\end{align*}
with some $c_k>0$.
Here in the last line,
 we used that $N_i=\lceil T_i s \rceil \ge s=\lceil \delta^{-CL}  \rceil$. We choose $C=C(k) > C_0 >0$, so that 
\begin{equation}\label{eq:c_0}
\delta^2+c_k\delta^{CL}\le 3\delta/4.
\end{equation}
Then since we are assuming that (\ref{eq:pos1}) fails, we deduce from the previous estimate that
\begin{equation}\label{eq:contr}
\left| \frac{1}{N_1\cdots N_k}\sum_{n\in\mathcal{N}} f(p(n)\Lambda)-\int_X f\,d\mu\right|
\ge (\delta-\delta^{L}-c_k\delta^{CL})\|f\|_{Lip}\ge \delta/4\|f\|_{Lip}.
\end{equation}

Now we apply \cite[Th.~8.6]{Green-Tao} to the polynomial map $p(n)$.
Note that
$$
\pi\left(p(n)\right)=D\pi\left(u+v+\sum_{i=1}^k\frac{n_i}{s}w_i\right)+\pi(g),
$$
and
$$
\pi\left(p(n)\right)-\pi\left(p(n-e_i)\right)=\frac{D\pi(w_i)}{s}.
$$

By \cite[Th.~8.6]{Green-Tao},
there exist $L_1,L_2>0$ such that for every $\rho\in (0,1/2)$ and $N_1,\ldots, N_k\ge 1$, one of the following holds:
\begin{enumerate}
\item[(i$'$)] For every Lipschitz function $f:X\to \mathbb{R}$,
\begin{equation}\label{eq:green-tao1}
\left|\frac{1}{N_1\cdots N_k}\sum_{n\in\mathcal{N}} f(p(n)\Lambda) - \int_X f\, d\mu\right|\le \rho\|f\|_{Lip}.
\end{equation}
\item[(ii$''$)] There exists $z\in \mathbb{Z}^l\backslash \{0\}$ such that
\begin{equation}\label{eq:green-tao2}
\|z\|\ll \rho^{-L_1}\quad\hbox{and}\quad \hbox{dist}\left(\left<z,\frac{D\pi(w_i)}{s}\right>,\Z\right)\ll
\rho^{-L_2}/N_i,\quad i=1,\ldots, k,
\end{equation}
where the implied constants depend only on the degree of the polynomial map.
\end{enumerate}
Comparing \eqref{eq:contr} and  \eqref{eq:green-tao1}, we deduce that
(ii$''$) holds with $\rho=\delta/4$, and there exists 
$z\in \mathbb{Z}^l\backslash \{0\}$ such that
\begin{equation}\label{eq:cond0}
\|z\|\ll \delta^{-L_1}\quad\hbox{and}\quad \hbox{dist}\left(\left<z,\frac{D\pi(w_i)}{s}\right>,\Z\right)\ll
\delta^{-L_2}/N_i,\quad i=1,\ldots,k.
\end{equation}
Since $\|w_i\|\le 1$, we obtain
\begin{equation}\label{eq:e1}
\left|\left<z,\frac{D\pi(w_i)}{s}\right>\right|\le \|z\|\|D\pi\| \|w_i\| s^{-1}\ll \delta^{-L_1+CL}\le \delta^{-L_1+C}.
\end{equation}
Taking $C=C(L_1)>0$ sufficiently large, the above estimate implies that
\begin{equation}\label{eq:e2}
\left|\left<z,\frac{D\pi(w_i)}{s}\right>\right|\le 1/4.
\end{equation}
Then
$$
\hbox{dist}\left(\left<z,\frac{D\pi(w_i)}{s}\right>,\Z\right)=\left|\left<z,\frac{D\pi(w_i)}{s}\right>\right|,
$$
and it follows from (\ref{eq:cond0}) that
$$
\left|\left<z,D\pi(w_i)\right>\right|\ll s\delta^{-L_2}/N_i\le \delta^{-L_2}/T_i,\quad i=1,\ldots,k.
$$
Hence, (\ref{eq:pos2}) holds, as required. This completes the proof of the theorem.
\QED

%We call the box as in (\ref{eq:box}) cube if the basis $\{v_i\}$ is orthogonal and $\|v_1\|=\cdots=\|v_k\|$.

We call a box map $\iota$, defined as in (\ref{eq:box}), {\it $(c_1,c_2)$-Diophantine} if 
there exists at least one  vector $w\in  \Omega := [-1,1]D\pi(w_1)+\cdots +[-1,1]D\pi (w_k)$ such that 
\begin{equation}\label{eq:diophh0}
\left|\left<z,w\right>\right|\ge c_1\|z\|^{-c_2}\quad \hbox{for all $z\in\Z^l\backslash\{0\}$.}
\end{equation}

We emphasize  that only one element of $\Omega$ has to satisfy the relevant Diophantine condition.  This
allows for the following remark which we will use later, e.g. in the proof of Theorem
\ref{th:mixing_lipschitz}.

\begin{remark}  \label{rotated box}
{\em
Let $\iota$ be a {\it $(c_1,c_2)$-Diophantine} box map,  $W$ the subspace spanned by the image of
$\iota$, and $S$ a compact subset of $\hbox{GL}(W)$.  Then there exists a constant $c = c(S)>0$,
which only depends
on  $S$, such that for all $s \in S$,  the box map $s \circ \iota$  is $(c \,c_1,c_2)$-Diophantine.  Indeed, 
since $S$ is compact, there exists $c = c(S)>0$ such that for every $s\in S$,
$$
[-1,1]D\pi(w_1)+\cdots +[-1,1]D\pi (w_k)\subset [-c^{-1},c^{-1}]D\pi(sw_1)+\cdots +[-c^{-1},c^{-1}]D\pi
(sw_k).
$$
If  $w\in [-1,1]D\pi(w_1)+\cdots +[-1,1]D\pi (w_k)$ satisfies
(\ref{eq:diophh0}), then $cw\in [-1,1]D\pi(s w_1)+\cdots +[-1,1]D\pi (s w_k)$
and satisfies (\ref{eq:diophh0}) with $c_1$ replaced by $c\,c_1$.
Hence, the box map $s \circ \iota$  is $(c \,c_1,c_2)$-Diophantine.
}
\end{remark}

The following corollary will play a crucial role in the next section.

\begin{corollary}\label{cor:line}
Given  $\theta,c_1,c_2>0$, there exists $\kappa=\kappa(c_2,\theta)>0$ such that
for every $\theta$-H\"older function $f:X\to \R$, 
$u\in \mathcal{L}(G)$, $(c_1,c_2)$-Diophantine box map $\iota:B\to\mathcal{L}(G)$, and $x\in X$,
we have
$$
\frac{1}{|B|}\int_B f(\exp(u)\exp(\iota(t))x)\, dt=\int_X f\,d\mu+O_{c_1,c_2}(\min(B)^{-\kappa}\|f\|_{C^\theta}).
$$ 
\end{corollary}

\proof
We first give a proof assuming that the function $f$ is Lipschitz.

We write the box map $\iota$ as
$$
\iota(t)=v+t_1 w_1+\cdots+ t_k w_k,\quad t\in B=[0,T_1]\times\cdots \times [0,T_k]
$$
with $v,w_1,\ldots,w_k\in W$ and $T_1,\ldots,T_k>0$.

We take $\kappa,\epsilon>0$ such that $\frac{-L_2\kappa+1}{L_1\kappa}> c_2$ and moreover
$\frac{-L_2(\kappa+\epsilon)+1}{L_1\kappa}> c_2$, where $L_1$ and $L_2$ are as in Theorem \ref{th:green-tao}.
Let $\delta=\min(B)^{-\kappa}$. 
We first assume that $\min(B)$ is sufficiently large, so that $\delta<1/2$.
Then by Theorem \ref{th:green-tao}, either 
\begin{equation}\label{eq:11}
\left| \frac{1}{|B|}\int_{B} f(\exp(u)\exp(\iota(t))x)\, dt-\int_X f\,d\mu\right|\le \min(B)^{-\kappa}\|f\|_{Lip}
\end{equation}
for all Lipschitz functions $f:X\to \R$, $u\in\mathcal{L}(G)$ and $x\in X$, or 
there exists $z\in \mathbb{Z}^l\backslash \{0\}$ such that
\begin{align*}\label{eq:22}
\|z\|&\ll \min(B)^{L_1\kappa},\\
|\left<z,D\pi(w_i)\right>|&\ll \min(B)^{L_2\kappa}/T_i\le \min(B)^{L_2\kappa-1},\quad i=1,\ldots,k.
\end{align*}
If the latter holds, then we deduce that there exists $z\in \mathbb{Z}^l\backslash \{0\}$ such that
\begin{align*}
|\left<z,D\pi(w_i)\right>|&\ll \min(B)^{-L_2\epsilon} \min(B)^{L_2(\kappa+\epsilon)-1}\ll \min(B)^{-L_2\epsilon}
\|z\|^{-\frac{-L_2(\kappa+\epsilon)+1}{L_1\kappa}}\\
&\le \min(B)^{-L_2\epsilon} \|z\|^{-c_2}
\end{align*}
for all $i=1,\ldots,k$. Writing $w=\sum_{i=1}^k a_i D\pi(w_i)$ with $a_i\in [-1,1]$, we also deduce that
$$
|\left<z,D\pi(w)\right>|\le
\sum_{i=1}^k  |\left<z,D\pi(w_i)\right>|
\ll \min(B)^{-L_2\epsilon} \|z\|^{-c_2}.
$$
When $\min(B)$ is sufficiently large, this estimate contradicts (\ref{eq:diophh0}).
Hence, we conclude that when $\min(B)\ge T_0=T_0(c_1,c_2)$, (\ref{eq:11}) holds and
\begin{align*}
\frac{1}{|B|}\int_{B} f(\exp(u)\exp(t)x)\, dt=\int_X f\,d\mu+O(\min(B)^{-\kappa}\|f\|_{Lip}).
\end{align*}
It is also clear that this estimate holds in the range $[0,T_0]$
with the implicit constant depending on $T_0$, and this completes proof of the corollary for Lipschitz functions.

For \holder functions,  we use the following well-known approximation result.  While we only use the
estimate of the Lipschitz norm here, we will  need  this lemma in full in Section~\ref{sec:cohomology}.

 \begin{Lemma} \label{convolution-approximation}
 Given  $\varepsilon >0$ and $0< \theta \leq 1$,  for any $\theta$-H\"{o}lder function $f:X   \to \R$,
 there is a $C^{\infty}$ function $f_{\varepsilon }:X\to\mathbb{R}$ which satisfies the following bounds 
   \begin{equation} \label{feg}
 \|f_{\varepsilon } -f \|_{C^0} \leq \varepsilon ^{\theta} \|f\|_{C^\theta}  \quad \text{and} \quad
   \|f_\epsilon\|_{Lip} \ll \epsilon^{-\dim(X)-1} \|f\|_{C^0}.  
   \end{equation}
     Furthermore,  for all $ l \in \mathbb N  $,
     \begin{equation} \label{feg-smooth}
\|f_{\varepsilon}\|_{C^l}  \ll_l \, \varepsilon ^{-\dim (X) -l}  \|f\| _{C^0}.
\end{equation}
 \end{Lemma}

 \proof  Given a $\theta$-H\"older function $f:X\to\mathbb{R}$, we set 
$$
f_\epsilon(x):=\int_G \phi_\epsilon(g^{-1}) f(gx)\, dm(g),
$$
where $m$ denotes the Haar measure on $G$, and $\phi_\epsilon$ is a nonnegative function such that
$$
\|\phi_\epsilon\|_{Lip}\ll \epsilon^{-\dim(X)-1},\quad \int_G\phi_\epsilon\, dm=1,\quad
\hbox{supp}(\phi_\epsilon)\subset B_\epsilon(e).
$$
Then
\begin{align*}
\|f_\epsilon-f\|_{C^0}\le \max_{x\in X} \int_G \phi_\epsilon(g^{-1})|f(gx)-f(x)|\, dm(g)\le \epsilon^\theta \|f\|_{C^\theta}.
\end{align*}
For $x,y\in X$ satisfying $d(x,y)<\epsilon$, we can write $y=hx$ with $h\in B_\epsilon(e)$. Then 
\begin{align*}
|f_\epsilon(x)-f_\epsilon(y)|\le \int_G |\phi_\epsilon(g^{-1})-\phi_\epsilon(hg^{-1})|f(gx)|\, dm(g)\ll
\epsilon^{-\dim(X)-1} \|f\|_{C^0}.
\end{align*}
Hence,
$$
\|f_\epsilon\|_{Lip} \ll \epsilon^{-\dim(X)-1} \|f\|_{C^0}.
$$
We can further assume that 
$\phi _{\varepsilon} $
 satisfies for all $ l \in \mathbb N $,
$$
\|  \phi _{\varepsilon} \| _{C^l}\ll_l \varepsilon ^{-\dim (X)-l}  \|\phi \| _{C^l},
$$
and it follows that 
 \begin{equation*} 
\|f_{\varepsilon}\|_  {C^l} \ll_l  \varepsilon ^{-\dim(X)-l}  \|f\| _{C^0}, 
\end{equation*}
as  the lemma claims.   \QED

Returning to the proof of Corollary \ref {cor:line},
 we obtain
\begin{align*}
\frac{1}{|B|}\int_{B} f(\exp(u)\exp(t)x)\, dt &=
\frac{1}{|B|}\int_{B} f_\epsilon(\exp(u)\exp(t)x)\, dt+O(\epsilon^\theta \|f\|_{C^\theta})\\
&=\int_X f_\epsilon\, d\mu+O\left(\min(B)^{-\kappa} \|f_\epsilon\|_{Lip}+\epsilon^\theta \|f\|_{C^\theta}\right)\\
&=\int_X f\, d\mu+O\left((\epsilon^{-\dim(X)-1} \min(B)^{-\kappa}+\epsilon^\theta) \|f\|_{C^\theta}\right).
\end{align*}
To optimise the error term, we set $\epsilon=\min(B)^{-\kappa/(\dim(X)+\theta+1)}$.
We readily obtain the corollary for H\"older functions.
\QED

We remark that the procedure just outlined applies quite generally, and allows  to go from estimates for Lipschitz functions to ones for \holder functions.  In particular, exponential mixing for Lipschitz or even only smooth functions always implies exponential mixing for \holder functions.

%%%%%%%%%%%%%%%%%%%%%%%%%%%%%%%%%%
%%%%%%%%%%%%   MIXING  %%%%%%%%%%%

\section{Mixing}\label{sec:mix}

In this section, we prove  Theorem \ref{exp mixing}    on  exponential mixing.  Let us recall the statement: 

\begin{Theorem}\label{th:mixing_lipschitz}
Let $\alpha$ be an ergodic automorphism of a compact nilmanifold $X=G/\Lambda$.
Then there exists $\rho=\rho(\theta)\in (0,1)$ such that for all 
$\theta$-H\"older functions $f_0,f_1: X \to \R$ and $n\in \mathbb{N}$,
\begin{align*}
 \int_X f_0(x)f_1(\alpha^n(x))\,d\mu(x)=\left(\int _X {f}_0\,  d \mu\right)\left(\int _X f_1\, d
   \mu\right)
+O\left(\rho ^n \|f_0 \|_{C^\theta}\|f_1\|_{C^\theta}  \right).
\end{align*}
\end{Theorem}

We denote by $\mu$ the Haar probability measure on $X$, and by $m$ the Haar measure on $G$
which is normalised, so that $m(F)=1$ where $F$ is a fundamental domain for $G/\Lambda$.

Every automorphism $\beta$ of $G$ defines a Lie-algebra automorphism $D\beta:\mathcal{L}(G)\to
\mathcal{L}(G)$ such that $\beta\circ \exp=\exp\circ D\beta$.
If $\beta(\Lambda)\subset \Lambda$, then $D\beta$ preserves the
rational structure of $\mathcal{L}(G)$ defined by $\Lambda$.

As in Section \ref{box maps}, we equip the group $G$ with the structure of an algebraic group, so that 
$\exp$ is a polynomial isomorphism.
More precisely, one can construct a basis, a so-called Malcev basis, $\{e_1,\ldots,e_d\}$ of 
$\mathcal{L}(G)_{\mathbb{Q}}$, such that the map
$$
\mathbb{R}^d\to G:(t_1,\ldots,t_d)\mapsto \exp(t_1e_1)\cdots\exp(t_de_d)
$$
is a polynomial isomorphism, 
$$
\Lambda=\exp(\mathbb{Z}e_1)\cdots\exp(\mathbb{Z}e_d),
$$
and
$$
F:=\exp([0,1)e_1)\cdots\exp([0,1)e_d)\subset G
$$
is a fundamental domain for $G/\Lambda$ (see \cite[1.2.7, 5.1.6, 5.3.1]{CG}).

We present the proof of Theorem \ref{th:mixing_lipschitz} in two stages: in Section \ref{sec:mix1}, we
give a proof assuming a suitable irreducibility condition, and in Section \ref{sec:mix2},
we reduce the proof to the irreducible case using an inductive argument.
 
\subsection{ Proof under an irreducibility assumption}\label{sec:mix1}
Let $w$ be a (real or complex) eigenvector of $D\alpha$ acting on  $\mathcal{L}(G)\otimes {\mathbb{C}}$ with eigenvalue
$\lambda$ such that $|\lambda|>1$. Such an eigenvector exists by the following lemma.

\begin{lemma}\label{l:eigenvalue}
If $\alpha$ is an ergodic automorphism of a nontrivial compact nilmanifold $X=G/\Lambda$,
then $D\alpha$ has an eigenvalue $\lambda$ with $|\lambda|>1$.
\end{lemma}

\proof
By \cite[5.4.13]{CG}, $\Lambda G'/G'$ is a lattice in $G/G'\simeq \mathbb{R}^l$.
The automorphism $\alpha$ defines a linear 
automorphism of the torus $T:=G/(\Lambda G')\simeq \mathbb{R}^l/L$,
where $L$ is a lattice in $\mathbb{R}^l$, and there is an $\alpha$-equivariant map $X\to T$ induced by $\pi$.
Since $\alpha|_{\mathbb{R}^l}$ preserves the lattice $L$, it follows that the eigenvalues
of $\alpha|_{\mathbb{R}^l}$ are algebraic integers. If we suppose that all these eigenvalues
satisfy $|\lambda|\le 1$, then it follows from \cite[Th.~1.31]{EW} that all the eigenvalues
of $\alpha|_{\mathbb{R}^l}$ are roots of unity. Then the automorphism
$\alpha|_T$ is not ergodic, and this contradicts
ergodicity of $\alpha$. Hence, $\alpha|_{\mathbb{R}^l}$ has an eigenvalue $\lambda$ with $|\lambda|>1$,
and this implies that $D\alpha$ has such an eigenvalue as well.
\QED

Since $D\alpha$ preserves the rational structure on $\mathcal{L}(G)$
defined by the lattice $\Lambda$, we may choose the eigenvector $w$ with coordinates in the algebraic closure 
$\overline{\mathbb{Q}}$.
In the real case, we denote by $W$ the corresponding one-dimensional eigenspace of $\mathcal{L}(G)$.
In the complex case, we denote by $W$ the two-dimensional subspace $\left<w,\bar w\right>\cap
\mathcal{L}(G)$, where $\bar w$ denotes the complex conjugate. We note that in a suitable basis
\begin{equation}\label{eq:isom}
D\alpha|_W=r\cdot \omega
\end{equation}
where $r=|\lambda|>1$ and $\omega$ is a rotation by angle $\hbox{Im}(\lambda)$.

In this subsection, we give a proof of Theorem \ref{th:mixing_lipschitz}
assuming that $D\pi(W)$ is not contained in any proper rational subspace of $\mathbb{R}^l$. 
This condition is used to guarantee existence of a ``generic'' vector in $D\pi(W)$ given by the following lemma.

\begin{lemma}\label{l:subspace}
Let $V\subset \mathbb{R}^l$ be a subspace defined over $\overline{\mathbb{Q}}\cap \mathbb{R}$ such that $V$
is not contained in any proper subspace defined over $\mathbb{Q}$.
Then there exists $w\in V\cap \overline{\mathbb{Q}}^l$ whose coordinates are real numbers   linearly independent over $\mathbb{Q}$.
\end{lemma}

\proof
Let $\{v_1,\ldots,v_s\}$ be a basis of $V$ whose coordinates $v_{ij}$ are in $\overline{\mathbb{Q}}\cap \mathbb{R}$.
We denote by $K$ the field generated by these coordinates. Clearly, $K$ is a finite extension of
$\mathbb{Q}$.
We can pick $\alpha_1,\ldots, \alpha_s\in \overline{\mathbb{Q}}\cap \mathbb{R}$ which are linearly
independent over $K$ (for instance, we can take a sufficiently large finite extension $K'$ of $K$
and choose $\{\alpha_i\}$ from a basis of $K'$ over $K$).

We set $w=\sum_{i=1}^s \alpha_i v_i$.
Suppose that there exists $c\in \mathbb{Q}^l$ such that $c\cdot w=0$. Then we have
$$
c\cdot w=\sum_{j=1}^l c_j\left(\sum_{i=1}^s\alpha_iv_{ij}\right)
=\sum_{i=1}^s\left(\sum_{j=1}^l c_j v_{ij}\right)\alpha_i=0.
$$
Now because $\sum_{j=1}^l c_j v_{ij}$ is in $K$, it follows that
$\sum_{j=1}^l c_j v_{ij}=0$ for all $i$, and $c\cdot V=0$. 
Since $V$ is not contained in any proper rational subspace, we conclude that $c=0$, which concludes the
proof.
\QED

As we remarked above, the subspace $W$ is defined over $\overline{\mathbb{Q}}$.
Moreover, since $W$ is invariant under complex conjugation, it is defined over 
$\overline{\mathbb{Q}}\cap \mathbb{R}$. This implies that 
the subspace $D\pi(W)$ is also defined over $\overline{\mathbb{Q}}\cap \mathbb{R}$.
Hence, by Lemma \ref{l:subspace}, $D\pi(W)$
contains a vector $w$ whose coordinates are real algebraic numbers
that are linearly independent over $\mathbb{Q}$.
By \cite[Th.~7.3.2]{BG}, there exist $c_1,c_2>0$ (in fact, one can take any $c_2>l-1$) such that
\begin{equation}\label{eq:DDD}
|\left<z,w \right>|\ge c_1 \|z\|^{-c_2}\quad\hbox{for all $z\in\mathbb{Z}^l\backslash \{0\}$,}
\end{equation}
This will allow us to apply Corollary \ref{cor:line} to box maps $\mathbb{R}^{\dim(W)}\to W$.

Let $E\subset\mathcal{L}(G)$ be the preimage of the fundamental domain $F$
under the exponential map. Since $E$ is the image of $[0,1)^d$ under a polynomial isomorphism,
it is a domain in $\mathcal{L}(G)$ with a piecewise smooth boundary.
We fix a basis of $\mathcal{L}(G)$ which contains the  basis of $W$
and consider a tessellation of $\mathcal{L}(G)$ by cubes $C$ of size $\epsilon$
with respect to this basis. Then
\begin{equation}\label{eq:region}
\left|E-\bigcup_{C\subset E} C \right|\ll \epsilon.
\end{equation}

Using the above notation, we rewrite the original integral as
\begin{align}\label{eq:est0}
\int_X f_0(x)f_1(\alpha^n(x))\, d\mu(x)&=
\int_F f_0(g\Lambda)f_1(\alpha^n(g)\Lambda)\, dm(g)\\
&=
\int_E f_0(\exp(u)\Lambda) f_1(\exp((D\alpha)^nu) \Lambda)\, du,\nonumber
\end{align}
where we used that the Haar measure on $G$ is the image of a suitably
normalised Lebesgue measure on $\mathcal{L}(G)$ under the exponential map \cite[1.2.10]{CG}.
It follows from (\ref{eq:region}) that
\begin{align}\label{eq:est1}
&\int_E f_0(\exp(u)\Lambda) f_1(\exp((D\alpha)^nu) \Lambda)\, du\\
=&\sum_{C\subset E} \int_C f_0(\exp(u)\Lambda) f_1(\exp((D\alpha)^nu) \Lambda)\, du
+O(\epsilon \|f_0\|_{C^0}\|f_1\|_{C^0}). \nonumber
\end{align}
For every cube $C$ in the above sum, we fix $u_C\in C$. Then for all $u\in C$,
$$
|f_0(\exp(u)\Lambda)-f_0(\exp(u_C)\Lambda)|\le d(\exp(u),\exp(u_C))\|f_0\|_{Lip}\ll \epsilon^\theta \|f_0\|_{C^\theta},
$$
and 
\begin{align}\label{eq:est2}
&\int_C f_0(\exp(u)\Lambda) f_1(\exp((D\alpha)^nu) \Lambda)\, du\\
=& f_0(\exp(u_C)\Lambda) \int_C f_1(\exp((D\alpha)^nu) \Lambda)\, du
+O(\epsilon^\theta \|f_0\|_{C^\theta}\|f_1\|_{C^\theta}). \nonumber
\end{align}
Since the cubes $C$ are chosen in a compatible way with the subspace $W$,
they can be written as $C=B'+B$ where $B$ is a cube in $W$
and $B'$ is a cube in the complementary subspace.
Given a cube $B\subset W$, we introduce a box map $\iota_B:\mathbb{R}^{\dim(W)}\to W$,
defined with respect to the fixed basis of $W$, such that $\iota_B([0,\epsilon]^{\dim(W)})=B$.
Since $\omega$ is a rotation, it follows from Remark \ref{rotated box} that for some $c>0$, each of the
box maps 
$$
\mathbb{R}^{\dim(W)}\to W:t\mapsto v+ \omega^n\iota_B(t),\quad v\in\mathcal{L}(G),
$$ 
is $(c\, c_1,c_2)$-Diophantine.
Therefore, applying Corollary \ref{cor:line}, we obtain there exists $\kappa>0$
such that for every $v\in \mathcal{L}(G)$, 
\begin{align}\label{eq:main:cor}
\frac{1}{|B|}\int_B f_1(\exp(v+(D\alpha)^nb) \Lambda)\, db
&=\epsilon^{-\dim(W)}\int_{[0,\epsilon]^{\dim(W)}} f_1(\exp(v+(D\alpha)^n \iota_B(t)) \Lambda)\, dt\\
&=(r^n\epsilon)^{-\dim(W)}\int_{[0,r^n\epsilon]^{\dim(W)}} f_1(\exp(v+\omega^n\iota_B(t)) \Lambda)\,
dt\nonumber \\
&=\int_X f_1\, d\mu
+O\left((r^n\epsilon)^{-\kappa}\|f_1\|_{C^\theta}\right). \nonumber
\end{align} 
Since this estimate is uniform over $v\in\mathcal{L}(G)$, we deduce that
\begin{align*}
\frac{1}{|C|}
\int_C f_1(\exp((D\alpha)^nu) \Lambda)\, du&=\frac{1}{|B'||B|}\int_{B'}\int_Bf_1(\exp((D\alpha)^nb'+ (D\alpha)^nb)
\Lambda)\, dbdb'\\
&=\int_X f_1\, d\mu+O\left((r^n\epsilon)^{-\kappa}\|f_1\|_{C^\theta}\right).
\end{align*}
Combining the last estimate with (\ref{eq:est1}) and (\ref{eq:est2}), we deduce that
\begin{align*}
\int_E f_0(\exp(u)\Lambda) f_1(\exp((D\alpha)^nu) \Lambda)\, du
=&\left(\sum_{C\subset E} f_0(\exp(u_C)\Lambda)|C|\right)\int_X f_1\, d\mu\\
&+O\left(\left(\sum_{C\subset E} |C| (r^n\epsilon)^{-\kappa} +\epsilon^\theta\right)\|f_0\|_{C^\theta}\|f_1\|_{C^\theta}\right).
\end{align*}
Since $f_0$ is $\theta$-H\"older and $\hbox{diam}(C)\ll \epsilon$, we obtain using (\ref{eq:region}),
\begin{align}\label{eq:lipschitz}
\sum_{C\subset E} f_0(\exp(u_C)\Lambda)|C|
&= \sum_{C\subset E} \int_C f_0(\exp(u)\Lambda)\, du+O(\epsilon^\theta \|f_0\|_{C^\theta})\\
&=\int_E f_0(\exp(u)\Lambda)\, du+O(\epsilon^\theta \|f_0\|_{C^\theta})\nonumber\\
&=  \int_X f_0\, d\mu +O(\epsilon^\theta \|f_0\|_{C^\theta}).\nonumber
\end{align}
Hence,
\begin{align*}
\int_E f_0(\exp(u)\Lambda) f_1(\exp((D\alpha)^nu) \Lambda)\, du
=&\left(\int_X f_1\, d\mu\right) \left(\int_X f_0\, d\mu\right)\\
&+O\left((r^n\epsilon)^{-\kappa}+\epsilon^\theta) \|f_0\|_{C^\theta}\|f_1\|_{C^\theta}\right).
\end{align*}
To optimise the error term, we choose $\epsilon=r^{-n\kappa/(\kappa+\theta)}$. Then
\begin{align*}
\int_X f_0(x)f_1(\alpha^n(x))\, d\mu(x)&=\int_E f_0(\exp(u)\Lambda) f_1(\exp((D\alpha)^nu) \Lambda)\, du\\
&=\left(\int_X f_0\, d\mu\right) \left(\int_X f_1\, d\mu\right)+O\left(\rho^n \|f_0\|_{C^\theta}\|f_1\|_{C^\theta}\right),
\end{align*}
where $\rho=r^{-\kappa\theta/(\kappa+\theta)}\in (0,1)$.
This proves Theorem \ref{th:mixing_lipschitz} under the irreducibility assumption.

We also observe that Corollary \ref{cor:line} implies the following
stronger version of  estimate (\ref{eq:main:cor}): for every $h\in G$,
automorphism $\beta$ of $G$ such that $\beta=id$ on $G/G'$, and $v\in \mathcal{L}(G)$,
\begin{align*}
\frac{1}{|B|}\int_B f_1(h\beta(\exp(v+(D\alpha)^nt)) \Lambda)\, dt=
\int_X f_1\, d\mu+O\left((r^n\epsilon)^{-\kappa}\|f_1\|_{C^\theta}\right).
\end{align*}
Indeed, using that $\beta\circ \exp=\exp\circ D\beta$, we obtain
\begin{align*}
&\frac{1}{|B|}\int_B f_1(h\beta(\exp(v+(D\alpha)^nt)) \Lambda)\, dt\\
=&
(r^n\epsilon)^{-\dim(W)}\int_{[0,r^n\epsilon]^{\dim(W)}} f_1(\exp((D\beta)v+(D\beta)\omega^n\iota_B(t))\Lambda)dt.
\end{align*}
Since $(D\pi)(D\beta)=D\pi$, the box maps
$$
t\mapsto (D\beta)v+(D\beta)\omega^n\iota_B(t)
$$
are also $(c\, c_1,c_2)$-Diophantine,  and the same estimate as in (\ref{eq:main:cor}) holds.
Therefore, the above argument implies that
\begin{align}\label{eq:inductive_step}
\int_X f_0(x)f_1(h\, \beta(\alpha^n(x)))\, d\mu(x)=
\left(\int_X f_0\, d\mu\right) \left(\int_X f_1\, d\mu\right)+O\left(\rho^n \|f_0\|_{C^\theta}\|f_1\|_{C^\theta}\right)
\end{align}
uniformly on $h\in G$ and automorphisms $\beta$ which preserve $\Lambda$ and act trivially on $G/G'$.

\begin{remark}\label{rem:affine}
{\em
Let $\sigma:X\to X$ be an affine diffeomorphism of a compact nilmanifold $X$.
Then $\sigma(x)=g_1\alpha(x)$ for $g_1\in G$ and an automorphism $\alpha$, and
$\sigma^n(x)=g_n\alpha^n(x)$ for $g_n\in G$. Since the estimate (\ref{eq:inductive_step})
is uniform over $h\in G$, it also holds for affine diffeomorphisms. 
This allows to extend the main results of this paper to affine diffeomorphisms. 
}
\end{remark}

\subsection{Proof of mixing in general}\label{sec:mix2}
We prove Theorem \ref{th:mixing_lipschitz} in general using induction on the dimension of 
the nilmanifold $X$.

Let $w\in \mathcal{L}(G)\otimes \mathbb{C}$
be an eigenvector of the automorphism $D\alpha$ with eigenvalue $\lambda$ of maximal modulus.
Since $\alpha$ is ergodic, $|\lambda|>1$ by Lemma \ref{l:eigenvalue}.
We set $W=\mathcal{L}(G)\cap \left< w, \bar w\right>$.
Since $D\alpha|_W$ has eigenvalues $\lambda$ and $\bar\lambda$, it follows either that
$D\alpha|_{[W,W]}$ must have eigenvalues of modulus $|\lambda|^2>|\lambda|$, or  $[W,W]=0$.  Hence
$\exp(W)$ is an abelian Lie subgroup of $G$.
By \cite[Ch. 3, Sec. 5]{S}, there exists a closed connected normal subgroup $M$ containing $\exp(W)$
such that $M/(M\cap \Lambda)$ is compact, and  for almost every $g\in G$, we have
$\overline{\exp(W) g \Lambda}=Mg\Lambda.$
Replacing the lattice $\Lambda$ by $g\Lambda g^{-1}$, we may assume without loss of generality that 
\begin{equation}\label{eq:dense}
\overline{\exp(W) \Lambda}=M\Lambda.
\end{equation}

\begin{lemma}\label{l:m_subgroup}
\begin{enumerate}
\item[(i)] The group $M$ is $\alpha$-invariant.
\item[(ii)] Denoting by $\pi:M\to M/M'$ the factor map,
$D\pi(W)$ is not contained in any proper rational subspace of $\mathcal{L}(M/M')$.
\item[(iii)] $[G,M]<M'$.
\end{enumerate}
\end{lemma}
 
\proof
We note that the group $M$ can be described
as the smallest closed connected normal subgroup containing $\exp(W)$ and intersecting $\Lambda$ in a
lattice (\cite[Ch. 3, Sec. 5]{S}). Equivalently, $M$ is the smallest closed connected subgroup
whose Lie algebra $\mathcal{L}(M)$ is an ideal in $\mathcal{L}(G)$ that contains $W$ and is defined over
$\mathbb{Q}$ with respect to the rational structure defined by $\Lambda$.
To show that $M$ is invariant under $\alpha$, we observe that
$$
\mathcal{L}(\alpha(M))=D\alpha(\mathcal{L}(M))
$$
also satisfies the above properties,
and so does
$$
\mathcal{L}(M\cap\alpha(M))=\mathcal{L}(M)\cap D\alpha(\mathcal{L}(M)).
$$
Therefore,  $\alpha(M)=M$ by minimality of $M$ proving (i).

To prove (ii), we consider the torus factor $M\Lambda/\Lambda\to T:=M\Lambda/(\Lambda M')$
induced by the map $\pi$. If $\pi(W)$ is contained in a proper rational subspace of $\mathcal{L}(M/M')$,
then the image of $D\pi(W)$ in $T$ is not dense, which contradicts (\ref{eq:dense}). This shows (ii).

Since the vector $w$ has coordinates in $\overline{\mathbb{Q}}$, so does the vector
$D\pi(w)$. For $\sigma\in \hbox{Gal}(\overline{\mathbb{Q}}/\mathbb{Q})$,
we denote by $D\pi(w)^\sigma$ its Galois conjugate. Then 
 $\left<D\pi(w)^\sigma:\sigma\in \hbox{Gal}(\overline{\mathbb{Q}}/\mathbb{Q}) \right>$ is a  rational subspace, 
contains $D\pi(W)$ and, hence, cannot be a proper subspace. 
This shows that $\hbox{Gal}(\overline{\mathbb{Q}}/\mathbb{Q})$ acts transitively 
on the eigenvalues of $D\alpha$ in $V:=\mathcal{L}(M/M')$. In particular,
it follows that $V$ does not contain any proper rational subspaces invariant under $D\alpha$.
Now we consider the adjoint action $\hbox{Ad}$ of $G$ on $V$.
Since $G$ is nilpotent, the set $V^G$ of $G$-fixed
points in $V$ is not trivial. Since $V^G$ is $(D\alpha)$-invariant and rational, we conclude that
$V^G=V$. This implies that every $g\in G$,
$$
(\hbox{Ad}(g)-id)(\mathcal{L}(M))\subset \mathcal{L}(M)',
$$
and the last claim of the lemma follows.
\QED

The nilmanifold $X=G/\Lambda$ fibers over 
the nilmanifold $Y=G/(M\Lambda)$ with fibers isomorphic to $Z=M\Lambda/\Lambda\simeq M/(M\cap\Lambda)$,
and we have the disintegration formula
\begin{equation}\label{eq:fiber}
\int_X f\, d\mu=\int_{Y} \int_Z f(yz)\, d\mu_Z(z)d\mu_Y(y),\quad f\in C(X),
\end{equation}
where $\mu_Y$ and $\mu_Z$ denote the normalised invariant measures on $Y$ and $Z$ respectively.
Since the groups $M$ and $\Lambda$ are $\alpha$-invariant, $\alpha$ defines transformations of $Y$ and
$Z$, and we obtain
\begin{align}\label{eq:d}
\int_X f_0(x) f_1(\alpha^n(x))\, d\mu(x)
&=\int_{Y} \left(\int_Z f_0(yz)f_1(\alpha^n(y)\alpha^n(z))\, d\mu_Z(z)\right)d\mu_Y(y)\\
&=\int_{F} \left(\int_Z f_0(gz) f_1(\alpha^n(g)\alpha^n(z))\, d\mu_Z(z)\right)dm_F(g),\nonumber
\end{align}
where $F\subset G$ is a bounded fundamental domain for $G/(M\Lambda)$, and
$m_F$ denotes the measure on $F$ induced by $\mu_Y$.

We claim that for some fixed $\rho\in (0,1)$ and every $g\in F$, 
\begin{align}\label{eq:induction}
\int_Z f_0(gz)f_1(\alpha^n(g)\alpha^n(z))\, d\mu_Z(z)
=& \left(\int_{Z} f_0(gz)\, d\mu_Z(z)\right)
\left(\int_{Z} f_1(\alpha^n(g)z)\, d\mu_Z(z)\right) \\
&+O(\rho^n \|f_0\|_{C^\theta} \|f_1\|_{C^\theta}) \nonumber
\end{align}
uniformly on $g\in F$.
To prove the claim above, we write 
$$
\alpha^n(g)=am\lambda\quad\hbox{with $a\in F$, $m\in M$, $\lambda\in \Lambda$.}
$$
Then
$$
\int_Z f_0(gz) f_1(\alpha^n(g)\alpha^n(z))\, d\mu_Z(z)
=\int_Z f_0(gz) f_1(am\beta(\alpha^n(z)))\, d\mu_Z(z),
$$ 
where $\beta$ denotes the transformation of $Z$ induced by the automorphism
$m\mapsto \lambda m\lambda^{-1}$,  $m\in M$.
We note %that it follows from Lemma \ref{l:m_subgroup} 
that $\beta$ acts trivially on $M/M'$ by Lemma \ref{l:m_subgroup}.
Let
$$
\phi_0(z):=f_0(gz)\quad\hbox{and}\quad \phi_1(z):=f_1(a z)\quad \hbox{with $z\in Z$}.
$$
Since $g, a\in F$, we have
$$
\|\phi_0\|_{C^\theta}\ll \|f_0\|_{C^\theta}\quad\hbox{and}\quad \|\phi_1\|_{C^\theta}\ll \|f_1\|_{C^\theta},
$$
and since $a(M\Lambda)=\alpha^n(g)(M\Lambda)$,
$$
\int_Z \phi_1\, d\mu_Z=\int_Z f_1(\alpha^n(g)z)\, d\mu_Z(z).
$$
Therefore, it follows from (\ref{eq:inductive_step}) that there exists $\rho\in (0,1)$ such that
\begin{align*}
&\int_Z \phi_0(z) \phi_1(m\beta(\alpha^n(z))))\, d\mu_Z(z)\\
=&\left(\int_Z \phi_0\, d\mu_Z\right)\left(\int_Z \phi_1\, d\mu_Z\right)
+O(\rho^n\|\phi_0\|_{C^\theta}\|\phi_1\|_{C^\theta})\\
=&\left(\int_Z f_1(gz)\, d\mu_Z(z)\right) \left(\int_Z f_0(\alpha^n(g)z)\, d\mu_Z(z)\right)
+O(\rho^n\|f_0\|_{C^\theta}\|f_1\|_{C^\theta})
\end{align*}
uniformly over $g,a\in F$, $m\in M$, and automorphisms $\beta$ of $Z$ which act trivially on $M/M'$.
This proves the claim (\ref{eq:induction}), and we conclude that
\begin{equation}\label{eq:last_ind}
\int_X f_0(x)f_1(\alpha^n(x))\, d\mu(x) 
=\int_{Y} \bar f_0(y) \bar f_1(\alpha^n(y)) \,d\mu_Y(y)+O(\rho^n\|f_0\|_{C^\theta}\|f_1\|_{C^\theta}),
\end{equation}
where the functions $\bar f_i:Y\to\mathbb{R}$ are defined by $y\mapsto \int_Z f_i(yz)\, d\mu_Z(z)$.
We note that
$$
\int_Y \bar f_i\,d\mu_Y=\int_X f_i\, d\mu. %\quad\hbox{and}\quad \|\bar f_i\|_{C^\theta}\le \|f_i\|_{C^\theta}.
$$
Since $\dim(Y)<\dim(X)$, Theorem \ref{th:mixing_lipschitz} follows from (\ref{eq:last_ind}) by induction on dimension.

%%%%%%%%%%%%%%%%%%%%%%%%%%%%%%%%%%
%%%%%%%%%%%%   MULTIPLE MIXING  %%%%%%%%%%%

\section{Multiple mixing}\label{sec:multiple}

In this section, we prove Theorem \ref{thm_mixing_lipschitz} on multiple exponential  mixing.  Let us recall the statement:

\begin{Theorem}\label{th:m_mixing_lipschitz}
Let $\alpha$ be an ergodic automorphism of a compact nilmanifolds $X=G/\Lambda$.
Then there exists $\rho=\rho(\theta)\in (0,1)$ such that for all 
$\theta$-H\"older function $f_0,\ldots,f_s: X \to \R$ and $n_0,\ldots,n_s\in \mathbb{N}$,
\begin{align*}
 \int_X \left(\prod_{i=0}^s f_i(\alpha^{n_i}(x))\right)\,d\mu(x)=
\prod_{i=0}^s \left(\int _X {f}_i\,  d \mu\right)
+O\left(\rho ^{\min_{i\ne j} |n_i-n_j|} \prod_{i=0}^s\|f_i \|_{C^\theta}  \right).
\end{align*}
\end{Theorem} 

We note that without loss of generality, we may assume that $n_0=0$ and $0<n_1<\cdots <n_s$.

As a preparation for the proof, we establish a result regarding equidistribution of 
images of box maps that generalises Corollary \ref{cor:line}.
We call a box map, defined as in (\ref{eq:box}), {\it $c_0$-bounded}  if $\|w_i\|\le c_0$
for all $i=1,\ldots,k$.

\begin{prop}\label{p:mmult}
Given $c_0,c_1,c_2,\theta>0$, there exists $\kappa=\kappa(c_2,\theta)>0$
such that for all $\theta$-H\"older functions
$f_1,\ldots,f_s:X\to\mathbb{R}$, $u_1,\ldots, u_s\in\mathcal{L}(G)$, automorphisms
$\beta_1,\ldots,\beta_s$ of $G$ such that $\beta_i=id$ on $G/G'$,
$0<r_1<\cdots< r_s$, $c_0$-bounded and $(c_1,c_2)$-Diophantine box maps $\iota_1,\ldots, \iota_s:B\to \mathcal{L}(G)$,
and $x_1,\ldots,x_s\in X$, we have 
\begin{align*}
\frac{1}{|B|}\int_B \left(\prod_{i=1}^s f_i(\exp(u_i)\beta_i(\exp(\iota_i(r_it)))x_i)\right)\, dt
=&\prod_{i=1}^s \left(\int_X f_i\,d\mu\right)\\
&+ O_{c_0,c_1,c_2}\left(\sigma(B,r_1,\ldots,r_s)^{-\kappa} \prod_{i=1}^s\|f_i\|_{C^\theta}\right),
\end{align*}
where $\sigma(B,r_1,\ldots,r_s)=\min\{\min(r_1B),r_sr_{s-1}^{-1},\ldots, r_2r_1^{-1} \}$.
\end{prop}

\proof
We first note that using the approximation argument as in the proof of
Corollary \ref{cor:line}, one can reduce the proof of the proposition to the case when
the functions are Lipschitz. Since this part is very similar to the proof of Corollary \ref{cor:line},
we omit details, and assume right away that the $f_i$'s are Lipschitz.

The proof involves applying Theorem  \ref{th:green-tao}   to the nilmanifold $X^{s}=G^s/\Lambda^s$.  
Let $L_1, L_2  >0$ be the constants from this theorem.
To simplify notation, we write $\sigma=\sigma(B,r_1,\ldots,r_s)$.
Let $\delta=\sigma^{-\kappa}$ where $\kappa>0$ is chosen so that 
$\frac{-\kappa(L_1+L_2)+1}{L_1\kappa}>c_2$ and moreover
$\frac{-(\kappa+\epsilon)(L_1+L_2)+1}{L_1\kappa}>c_2$ for some fixed $\epsilon>0$.
First, we assume that $\sigma$ is sufficiently large so that $\delta\in (0,1/2)$.

We write the box maps $\iota_i$ as 
$$
\iota_i(t)=v_i+t_1w_i^{(1)}+\cdots+t_kw_i^{(k)},\quad t\in B=[0,T_1]\times\cdots \times [0,T_k],
$$
with $v_i,w_i^{(j)}\in \mathcal{L}(G)$ and $T_1,\ldots T_k,>0$ and
set 
\begin{align*}
&f=f_1\otimes\cdots\otimes f_s:X^s\to\mathbb{R},\\
& u=(u_1,\ldots,u_s)\in \mathcal{L}(G)^s,\\
& \iota:B\to \mathcal{L}(G)^s:t\mapsto (D\beta_1\iota_1(r_1 t),\ldots, D\beta_s\iota_s(r_s t)),\\
&x=(x_1,\ldots,x_s)\in X^s.
\end{align*}
Then
$$
\int_B \left(\prod_{i=1}^s f_i(\exp(u_i)\beta_i(\exp(\iota_i(r_it))))x_i\right)\, db
=\int_B f(\exp(u)\exp(\iota(t))x)\, dt.
$$
Applying Theorem \ref{th:green-tao}, we deduce that for every $\delta\in (0,1/2)$, either
\begin{equation}\label{eq:good}
\left|\frac{1}{|B|}\int_B f(\exp(u)\exp(\iota(t)))\, dt-\int_X f\, d\mu\right|
\le \delta \|f\|_{Lip},
\end{equation}
or there exists $(z_1,\ldots,z_s)\in (\mathbb{Z}^l)^s\backslash \{0\}$ such that 
\begin{equation}\label{eq:bad1}
\|z_1\|,\ldots,\|z_s\|\ll \delta^{-L_1}=\sigma^{\kappa L_1}
\end{equation}
and 
\begin{equation}\label{eq:bad2}
\left| \sum_{j=1}^s r_j \left<z_j, D\pi D\beta_j (w_j^{(i)})\right>\right|\ll \delta^{-L_2}/\min(B)
=\sigma^{\kappa L_2}/\min(B)
\quad
\hbox{for all $i=1,\ldots,k$.}
\end{equation}
We note that since $\beta_j=id$ on $G/G'$, we have $D\pi D\beta_j =D\pi$.

Suppose that (\ref{eq:bad1})--(\ref{eq:bad2}) holds.
Since $\|w_j^{(i)}\|\le c_0$ by  assumption, using the triangle inequality we deduce that
$$
\left| \left<z_s, D\pi(w_s^{(i)})\right>\right|\ll \sigma^{\kappa L_2}/\min(r_s B)
+\sum_{j=1}^{s-1}\sigma^{\kappa L_1}/(r_s r_{j}^{-1})\le \sigma^{\kappa(L_1+L_2)-1}
$$
for all $i=1,\ldots,k$.
Then by (\ref{eq:bad1}),
\begin{align*}
\left| \left<z_s, D\pi(w_s^{(i)})\right>\right|
&\ll \sigma^{-(L_1+L_2)\epsilon} \sigma^{(\kappa+\epsilon)(L_1+L_2)-1}
\ll \sigma^{-(L_1+L_2)\epsilon} \|z\|^{-\frac{-(\kappa+\epsilon)(L_1+L_2)+1}{L_1\kappa}}\\
&\le \sigma^{-(L_1+L_2)\epsilon} \|z\|^{-c_2}.
\end{align*}
Since the box map $\iota_s$ is $(c_1,c_2)$-Diophantine, there exists $w_s\in \sum_{i=1}^k [-1,1]D\pi(w_s^{(i)})$ which satisfies
(\ref{eq:diophh0}). On the other hand, it follows from the previous estimate that
$$
\left| \left<z_s, D\pi(w_s)\right>\right|\le \sum_{i=1}^k
\left| \left<z_s, D\pi(w_s^{(i)})\right>\right|\ll \sigma^{-(L_1+L_2)\epsilon} \|z\|^{-c_2}.
$$
When $\sigma$ is sufficiently large, this estimate contradicts (\ref{eq:diophh0}),
unless $z_s=0$. Hence, we deduce that $z_s=0$. 

Now we repeat the above argument and deduce from (\ref{eq:bad1})--(\ref{eq:bad2})
that
$$
\left| \left<z_{s-1}, D\pi(w_{s-1}^{(i)})\right>\right|\ll \sigma^{\kappa L_2}/\min(r_{s-1} B)
+\sum_{j=1}^{s-2} \sigma^{\kappa L_1}/(r_{s-1} r_{j}^{-1})\le \sigma^{\kappa(L_1+L_2)-1}
$$
for all $i=1,\ldots,k$,  and ultimately that $z_{s-1}=0$, when $\sigma$ is sufficiently large.
Hence, we conclude that $(z_1,\ldots,z_s)=0$ when $\sigma\ge \sigma_0=\sigma_0(c_0,c_1,c_2)$.
Therefore, in this range (\ref{eq:good}) holds with $\delta=\sigma^{-\kappa}$.
 This proves the claim of the proposition for sufficiently
large $\sigma$. It is also clear that this estimate holds in the range $[0,\sigma_0]$
with the implicit constant depending on $\sigma_0$. This completes the proof of the proposition.
\QED

\subsection{Multiple mixing under irreducibility assumption}
In this section, we prove Theorem \ref{th:m_mixing_lipschitz} under the irreducibility 
condition as in Section \ref{sec:mix1}. Namely, $W$ denotes a $(D\alpha)$-invariant
subspace of $\mathcal{L}(G)$ such that $D\pi(W)$ is not contained in a proper rational subspace
and (\ref{eq:isom}) holds.

As in (\ref{eq:est0}), we obtain
\begin{align*}
\int_X f_0(x)\left(\prod_{i=1}^s f_i(\alpha^{n_i}(x))\right)\, d\mu(x)
=
\int_E f_0(\exp(u)\Lambda)\left( \prod_{i=1}^s f_i(\exp((D\alpha)^{n_i}u)\Lambda)\right)\, du.
\end{align*}
As in Section \ref{sec:mix1}, we tessellate the region $E$ by cubes $C$ of size $\epsilon$
which are compatible with the subspace $W$ and get 
\begin{align}\label{eq:disintegrate}
&\int_E f_0(\exp(u)\Lambda)\left(\prod_{i=1}^s f_i(\exp((D\alpha)^{n_i}u)\Lambda)\right)\, du\\
=& \sum_{C\subset E} f_0(\exp(u_C)\Lambda)
\int_C \left(\prod_{i=1}^s f_i(\exp((D\alpha)^{n_i}u)\Lambda)\right)\, du +O\left(\epsilon^\theta \prod_{i=0}^s\|f_i\|_{C^\theta}\right),\nonumber
\end{align}
where $u_C\in C$.
Each cube $C$ can be written as $C=B'+B$ where $B$ is a cube in $W$ and $B'$ is a cube in the
complementary subspace. For every cube $B$, we take a box map 
$\iota_B:\mathbb{R}^{\dim(W)}\to W$ such that $\iota_B([0,\epsilon]^{\dim(W)})=B$.
Because $\omega$ is a rotation, there exists $c_0>0$ such that each of the box maps
$$
\mathbb{R}^{\dim(W)}\to W: t\mapsto v+\omega^n\iota_B(t),\quad v\in \mathcal{L}(G),\; n\in\mathbb{N},
$$
is $c_0$-bounded.  It was also observed in Section \ref{sec:mix1} that each of these maps is 
$(c_1,c_2)$-Diophantine. 
Hence, Proposition \ref{p:mmult} implies that there exists $\kappa\in (0,1)$ such that 
uniformly on $v_1,\ldots,v_s\in \mathcal{L}(G)$,
\begin{align}\label{eq:CC}
&\frac{1}{|B|}\int_B \left(\prod_{i=1}^s f_i(\exp(v_i+(D\alpha)^{n_i}b)\Lambda)\right)\,db\\
=&\epsilon^{-\dim(W)}\int_{[0,\epsilon]^{\dim(W)}}
\left(\prod_{i=1}^s f_i(\exp(v_i+r^{n_i}\omega^{n_i}\iota_B(t))\Lambda)\right)\,dt \nonumber\\
=&\prod_{i=1}^s \left(\int_X f_i\,d\mu\right)+O\left(\sigma^{-\kappa} \prod_{i=1}^s\|f_i\|_{C^\theta}\right),
\nonumber
\end{align}
where $\sigma=\min\{\epsilon r^{n_1}, r^{n_2-n_1},\ldots, r^{n_s-n_{s-1}}\}$.
Since this estimate is uniform over $v_i$'s, we conclude that
\begin{align*}
&\frac{1}{|C|}\int_C \left(\prod_{i=1}^s f_i(\exp((D\alpha)^{n_i}u)\Lambda)\right)\, du\\
=&\frac{1}{|B'||B|}\int_{B'}\int_B  \left(\prod_{i=1}^s f_i(\exp((D\alpha)^{n_i}b' +(D\alpha)^{n_i}b)\Lambda)\right)
\, dbdb'\\
=&\prod_{i=1}^s \left(\int_X f_i\,d\mu\right)+O\left(\sigma^{-\kappa} \prod_{i=1}^s\|f_i\|_{C^\theta}\right).
\end{align*}
Now it follows from (\ref{eq:disintegrate}) that 
\begin{align*}
&\int_E f_0(\exp(u)\Lambda)\left(\prod_{i=1}^s f_i(\exp((D\alpha)^{n_i}u)\Lambda)\right)\, du\\
=& \left(\sum_{C\subset E} f_0(\exp(u_C)\Lambda)|C|\right)\prod_{i=1}^s \left(\int_X f_i\,d\mu\right)
+O\left((\sigma^{-\kappa}+\epsilon^\theta) \prod_{i=0}^s\|f_i\|_{C^\theta}\right),
\end{align*}
and by (\ref{eq:lipschitz}), 
\begin{align*}
\int_E f_0(\exp(u)\Lambda)\left(\prod_{i=1}^sf_i(\exp((D\alpha)^{n_i}u)\Lambda)\right)\,du
=&\prod_{i=0}^s\left(\int_X f_i\,d\mu\right)\\
&+O\left((\sigma^{-\kappa}+\epsilon^\theta)\prod_{i=0}^s \|f_i\|_{C^\theta}\right).
\end{align*}
Finally, taking $\epsilon=r^{-\kappa n_1/(\theta+\kappa)}$, we obtain
\begin{align*}
&\int_X f_0(x)\left(\prod_{i=1}^s f_i(\alpha^{n_i}(x))\right)\, d\mu(x)\\
=& \int_E f_0(\exp(u)\Lambda)\left(\prod_{i=1}^s f_i(\exp((D\alpha)^{n_i}u)\Lambda)\right)\, du\\
=& \prod_{i=0}^s \left(\int_X f_i\,d\mu\right)
+O\left(\min\{r^{\theta n_1/(\theta+\kappa)},r^{n_2-n_1},\ldots,r^{n_s-n_{s-1}}\}^{-\kappa} \prod_{i=1}^s\|f_i\|_{C^\theta}\right).
\end{align*}
This completes the proof of Theorem \ref{th:m_mixing_lipschitz} under the irreducibility assumption.

The proof of the general case will be given in the following section using an inductive argument.
For this purpose, we note that the above argument gives the following stronger result:
there exists $\rho\in (0,1)$ such that for every $h_1,\ldots,h_1\in G$ and automorphisms 
$\beta_1,\ldots,\beta_s$ of $G$ which preserve $\Lambda$ that act trivially on $G/G'$, we have 
\begin{align}\label{eq:step00}
&\int_X f_0(x)\left(\prod_{i=1}^s f_i(h_i\beta_i(\alpha^{n_i}(x)))\right)\, d\mu(x)\\
=& \prod_{i=0}^s\left(\int_X f_i\,d\mu\right)+O\left(\rho^{\min\{n_1,n_{2}-n_1,\ldots, n_s-n_{s-1}\}} 
\prod_{i=0}^s \|f_i\|_{C^\theta}\right)\nonumber
\end{align}
uniformly over $h_i$'s and $\beta_i$'s.
Indeed, Proposition \ref{p:mmult} implies that in (\ref{eq:CC}) we have, more generally,
\begin{align*}
&\frac{1}{|B|}\int_B \left(\prod_{i=1}^s f_i(h_i\beta_i(\exp(v_i+(D\alpha)^{n_i}b))\Lambda)
\right)\, db\\
=&\prod_{i=1}^s \left(\int_X f_i\,d\mu\right)+O\left(\sigma^{-\kappa} \prod_{i=1}^s\|f_i\|_{C^\theta}\right),
\end{align*}
and the rest of the proof can be carried out as well.

\subsection{Proof of multiple mixing in general}
We use notation introduced in Section \ref{sec:mix2}.
In particular, $W$ denotes a $(D\alpha)$-invariant subspace
of $\mathcal{L}(G)$, and we arrange that $\overline{\exp(W)\Lambda}=M\Lambda$ where 
$M$ is closed connected normal $\alpha$-invariant subgroup containing $\exp(W)$ such that
$M/(M\cap\Lambda)$ is compact. 

The nilmanifold $X=G/\Lambda$ fibers in $\alpha$-invariant fashion over 
the nilmanifold $Y=G/(M\Lambda)$ with fibers isomorphic to $Z=M\Lambda/\Lambda\simeq M/(M\cap\Lambda)$,
and the disintegration formula (\ref{eq:fiber}) holds.
Using this disintegration formula, we obtain, similarly to (\ref{eq:d}),
\begin{align}\label{eq:dis}
&\int_X f_0(x) \left(\prod_{i=1}^s f_i(\alpha^{n_i}(x))\right)\, d\mu(x)\\
=&\int_{Y} \left(\int_Z f_0(yz) \left(\prod_{i=1}^s f_i(\alpha^{n_i}(y)\alpha^{n_i}(z))\right)\,
  d\mu_Z(z)\right)d\mu_Y(y) \nonumber \\
=&\int_{F} \left(\int_Z f_0(gz)\left(\prod_{i=1}^s f_i(\alpha^{n_i}(g)\alpha^{n_i}(z))\right)\, d\mu_Z(z)\right)dm_F(g).\nonumber
\end{align}

We claim that there exists $\rho\in (0,1)$ such that for every $g\in F$, 
\begin{align}\label{eq:induction2}
&\int_Z f_0(gz)\left(\prod_{i=1}^s f_i(\alpha^{n_i}(g)\alpha^{n_1}(z))\right)\, d\mu_Z(z)\\
=& \left(\int_{Z} f_0(gz)\, d\mu_Z(z)\right)
\prod_{i=1}^s \left(\int_{Z} f_i(\alpha^{n_i}(g)z)\, d\mu_Z(z)\right)
+O\left(\rho^n \prod_{i=0}^s \|f_i\|_{C^\theta}\right) \nonumber
\end{align}
uniformly on $g\in F$.
To prove this claim, we write $\alpha^{n_i}(g)=a_{i}m_{i}\lambda_{i}$ with $a_{i}\in F$,
$m_{i}\in M$, and $\lambda_{i}\in \Lambda$. Then
\begin{align*}
\int_Z f_0(gz)\left(\prod_{i=1}^s f_i(\alpha^{n_i}(g)\alpha^{n_i}(z))\right)\, d\mu_Z(z)
= \int_Z f_0(gz) \left(\prod_{i=1}^s f_i(a_{i}m_{i}\beta_{i}(\alpha^{n_i}(z)))\right)\, d\mu_Z(z),
\end{align*}
where $\beta_{i}$ denotes the transformation of $Z$ induced by the automorphism
$m\mapsto \lambda_{i}m\lambda_{i}^{-1}$,  $m\in M$.
Note that by Lemma \ref{l:m_subgroup} the automorphism $\beta_{i}$ is trivial on $M/M'$.
Let
$$
\phi_0(z):=f_0(g z)\quad\hbox{and}\quad \phi_i(z):=f_i(a_i z),\; i=1,\ldots,s,\quad \hbox{with $z\in Z$}.
$$
Since $g$ and $a_i$'s belong to the compact set $F$, 
$$
\|\phi_i\|_{C^\theta}\ll \|f_i\|_{C^\theta},\quad i=0,\ldots s,
$$
and since $a_{i}(M\Lambda)=\alpha^{n_i}(g)(M\Lambda)$,
$$
\int_Z \phi_i\, d\mu_Z=\int_Z f_i(\alpha^{n_i}(g)z)\, d\mu_Z(z),\quad i=1,\ldots,s.
$$
Applying the estimate (\ref{eq:step00}), we deduce that for some $\rho\in (0,1)$,
\begin{align*}
&\int_Z \phi_0(z) \left(\prod_{i=1}^s \phi_i(m_{i}\beta_{i}(\alpha^{n_i}(z)))\right)\, d\mu_Z(z)\\
=&\prod_{i=1}^s\left(\int_Z \phi_i\, d\mu_Z\right)+O\left(\rho^n\prod_{i=0}^s\|\phi_i\|_{C^\theta}\right)\\
=&\left(\int_Z f_0(gz)\, d\mu_Z(z)\right) \prod_{i=1}^s \left(\int_Z f_i(\alpha^{n_i}(g)z)\, d\mu_Z(z)\right)
+O\left(\rho^n\prod_{i=0}^s \|f_i\|_{C^\theta}\right).
\end{align*}
This implies the claim (\ref{eq:induction2}).
Now combining (\ref{eq:induction2}) with (\ref{eq:dis}), we deduce that
\begin{align}\label{eq:last_ind0}
\int_X f_0(x)\left(\prod_{i=1}^s f_i(\alpha^{n_i}(x))\right)\, d\mu(x)
=&\int_{Y} \bar f_0(y) \left(\prod_{i=1}^s \bar f_i(\alpha^{n_i}(y))\right) \,d\mu_Y(y)\\
&+O\left(\rho^n\prod_{i=1}^s \|f_i\|_{C^\theta}\right).\nonumber
\end{align}
where the functions $\bar f_i:Y\to \mathbb{R}$ are defined by $y\mapsto \int_Z f_i(yz)\, d\mu_Z(z)$. Clearly,
$$
\int_Y \bar f_i\,d\mu_Y=\int_X f_i\, d\mu. %quad\hbox{and}\quad \|\bar f_i\|_{C^\theta}\le \|f_i\|_{C^\theta}.
$$
Since $\dim(Y)<\dim(X)$, Theorem \ref{th:m_mixing_lipschitz} now follows from (\ref{eq:last_ind0}) by induction on dimension.

%%%%%%%%%%%%%%%%%%%%%%%%%%%%%%%%%%
%%%%%%%%%%%%   EQUIDISTRIBUTION ON UNSTABLE MANIFOLDS  %%%%%%%%%%%

\section{Equidistribution of unstable manifolds}    \label{sec-unstable}

In this section we prove an equidistribution result for  unstable manifolds.  Besides its own intrinsic
interest, we will use this later in our treatment  of probabilistic limit theorems
in Section \ref{CLT}.

Let $\alpha$ be an ergodic automorphism of a compact nilmanifold $X=G/\Lambda$.
We denote by ${W}^\alpha\subset \mathcal{L}(G)$ the unstable subspace 
of $D\alpha$, namely, the subspace of $\mathcal{L}(G)$ spanned by Jordan subspaces
of $D\alpha$ with eigenvalues $\lambda$ satisfying $|\lambda|>1$.
Note that since $[{W}^\alpha,{W}^\alpha]\subset {W}^\alpha$, 
$\exp(W^\alpha)$ is a Lie subgroup of $G$. We decompose ${W}^\alpha$ as a direct sum
${W}^\alpha=\oplus_{i=1}^\ell {W}^\alpha_{i}$, so that 
$D\alpha|_{{W}^\alpha_{i}}$ acts as a (real) Jordan block. Namely, each subspace ${W}^\alpha_{i}$
has a basis $\{w_1,\ldots,w_s\}$ such that 
\begin{align}\label{eq:jordan1}
(D\alpha) w_i&=\lambda w_i+ w_{i+1},\,\,\,\, i<s,\\
(D\alpha) w_s&=\lambda w_s,\nonumber
\end{align}
where $\lambda$ is a real eigenvalue of $D\alpha$, or a basis 
$\{w_1,w_1',\ldots,w_s,w_s'\}$ such that 
\begin{align}\label{eq:jordan2}
(D\alpha) w_i&=a w_i+b w_i'+ w_{i+1},\,\,\, (D\alpha) w'_i=-b w_i+a w_i'+ w'_{i+1},\,\,\,\, i<s,\\
(D\alpha) w_s&=a w_s+b w_s',\,\,\, (D\alpha) w'_s=-b w_s+a w_s',
\end{align}
where $\lambda=a+bi$ is a complex eigenvalue of $D\alpha$.
We order the subspaces ${W}^\alpha_{i}$ with respect to the size of $|\lambda|$.
Then 
\begin{equation}\label{eq:nesting}
[{W}^\alpha,{W}_i^\alpha]\subset \oplus_{j>i}{W}_j^\alpha. 
\end{equation}
For each $i$, we define a map $\psi_i:\mathbb{R}^{\dim({W}_i^\alpha)}\to \exp(W_\alpha)$
which is either
$$
\psi_i: (t_1,\ldots,t_s)\mapsto \exp(t_1w_1)\cdots \exp(t_sw_s)
$$
in the real case, or 
$$
\psi_i: (t_1,t_1'\ldots,t_s,t_s')\mapsto \exp(t_1w_1+t'_1w'_1)\cdots \exp(t_sw_s+t_s'w_s')
$$
in the complex case.
Let $\psi:\mathbb{R}^{\dim({W}^\alpha)}\to \exp({W}^\alpha)$
be the product of the maps $\psi_i$. It follows from (\ref{eq:nesting})
that $\psi$ is a diffeomorphism  and that the image of the Lebesgue
measure gives the Haar measure on $\exp(W^\alpha)$ \cite[1.2.7]{CG}. 

\begin{Theorem}\label{th:unstable}
Let $\alpha$ be an ergodic automorphism of a compact nilmanifolds $X=G/\Lambda$.
Then there exist $\kappa=\kappa(\theta)>0$ and $\rho=\rho(\theta)\in (0,1)$ such that for every box $B\subset
\mathbb{R}^{\dim({W}^\alpha)}$,  $\theta$-H\"older function $f:X\to \mathbb{R}$,
$h\in G$, and  $g\in G$, we have
$$
\frac{1}{|B|}\int_B f(\alpha^n(h\psi(b))g\Lambda)\,db=\int_X f\, d\mu+O(\min(B)^{-\kappa}\rho^n \|f\|_{C^\theta}).
$$
\end{Theorem}

\proof
We give a proof using an  inductive argument similar to the proof of exponential mixing 
in Section \ref{sec:mix}.

Let $W={W}^\alpha\cap \left<w,\bar w\right>$ where $w$ is the eigenvector of $D\alpha$
in ${W}^\alpha_\ell$.
More explicitly, $W=\left<w_s\right>$ or $W=\left<w_s,w_s'\right>$ with notation
(\ref{eq:jordan1})--(\ref{eq:jordan2}).
As in Section \ref{sec:mix2}, we deduce that there exists a closed normal subgroup $M$ of $G$
containing $\exp(W)$ such that $M/(M\cap \Lambda)$ is compact and for almost all $g\in G$,
\begin{equation}\label{eq:dense0}
\overline{\exp(W)g\Lambda}=Mg\Lambda. 
\end{equation}

The map $\psi:\mathbb{R}^{\dim({W}^\alpha)}\to \exp({W}^\alpha)$ can be written as 
a product $\psi=\xi\cdot \eta$ with $\xi:\mathbb{R}^{\dim({W}^\alpha)-\dim(W)}\to \exp({W}^\alpha)$
and $\eta:\mathbb{R}^{\dim(W)}\to \exp(W)$, where 
$\eta:t\mapsto \exp(tw_s)$ or $\eta:(t,t')\mapsto \exp(tw_s+t'w_s')$
and $\xi$ is the product of the remaining exponential maps appearing in $\psi$. Then
$$
\int_B f(\alpha^n(h\psi(b))g\Lambda)\,db=\int_{C}\int_{D}
f(\alpha^n(h\xi(u)\eta(v))g\Lambda)\,dudv,
$$
where $C$ is a box in $\mathbb{R}^{\dim({W}^\alpha)-\dim(W)}$
and $D$ is a box in $\mathbb{R}^{\dim(W)}$ such that $B=C\times D$.

We first show that images of the map $\eta$ are equidistributed in a suitable sense.
Namely, we claim that there exists $\rho\in (0,1)$ such that for every $h\in G$ and every $g\Lambda\in X$
such that (\ref{eq:dense0}) holds,
\begin{equation}\label{eq:step}
\frac{1}{|D|}\int_{D} f(\alpha^n(h\eta(t))g\Lambda)\,dt=\int_{Z} f(\alpha^n(h)gm\Lambda)\, \mu_Z(m)
+O(\rho^n \|f\|_{C^\theta}),
\end{equation}
where $\mu_Z$ denotes the invariant normalised measure on the nilmanifold $Z=M/(M\cap \Lambda)$.
Let $F_0\subset G$ be a bounded subset such that $G=F_0\Lambda$.
Then there exists a bounded subset $F$ of $G$ such that $G=FM(g_0\Lambda g_0^{-1})$ for all $g_0\in F_0$.
Indeed, we can take $F=F_0F_0^{-1}$.  We note that in (\ref{eq:step}) we may assume that $g\in F_0$,
and to simplify notation, we replace $\Lambda$ by $g\Lambda g^{-1}$.
Then (\ref{eq:step}) holds with $g=e$.

Next we write
$\alpha^n(h)=am\lambda$ with 
$a\in F$, $m\in M$, and $\lambda\in\Lambda$. Then 
\begin{equation}\label{eq:step2}
\int_{D} f(\alpha^n(h\eta(t))\Lambda)\,dt
=\int_{D} f(am \beta(\alpha^n(\eta(t)))\Lambda)\,dt,
\end{equation}
where $\beta$ denotes the automorphism of $M$ defined by $m\mapsto
\lambda m\lambda^{-1}$. We note that $\beta$ acts trivially on $M/M'$ by Lemma \ref{l:m_subgroup}.
To analyse (\ref{eq:step2}), we apply Corollary \ref{cor:line} to the nilmanifold $Z=M/(M\cap \Lambda)$.
Setting $\phi(z):=f(az)$, $z\in Z$, we get 
$$
\int_{D} f(am \beta(\alpha^n(\eta(t)))\Lambda)\,dt=
\int_{D} \phi(m \beta(\alpha^n(\eta(t)))\Lambda)\,dt,
$$
and since $a\in F$, we have 
$$
\|\phi\|_{C^\theta}\ll \|f\|_{C^\theta}.
$$ 
For the next computation, let us assume that $\dim(W)=2$. When $\dim(W)=1$, the proof is similar and simpler.
We observe that $D\alpha|_W=r \omega$ where $r>1$ and $\omega$ is a rotation of $W$, so that
$$
\beta(\alpha^n(\eta(t,t')))=\exp(r^n(t(D\beta\omega^nw_s)+t'(D\beta\omega^nw'_s))).
$$
Making a change of variables, 
$$
\int_{D} \phi(m \beta(\alpha^n(\eta(t)))\Lambda)\,dt
= r^{-2n}\int_{r^n D} \phi(m \exp(\iota_n(t))\Lambda)\,dt,
$$
where $\iota_n$ denotes the box map $(t,t')\mapsto t(D\beta\omega^nw_s)+t'(D\beta\omega^nw'_s)$.
We note that $(D\pi)(D\beta)=D\pi$.
Since $\overline{\exp(W)\Lambda}=M\Lambda$, it follows that $D\pi(W)$ is not contained
in any proper rational subspace. In particular, it follows from Lemma \ref{l:subspace} and
\cite[Th.~7.3.2]{BG} that $D\pi(W)$ contains a vector
$w$ satisfying the Diophantine condition (\ref{eq:diophh0}). Since $\omega$ is an isometry, 
this implies that the box map $\iota_n$ is $(c_1,c_2)$-Diophantine where $c_1,c_2$ are uniform in $n$ and
$\beta$ (see Remark \ref{rotated box}). Therefore, Corollary \ref{cor:line} implies that there exists $\kappa>0$ such that
\begin{align*}
\frac{1}{|r^{n}D|}\int_{r^n D} \phi(m \beta(\exp(\iota_n(t)))\Lambda)\,dt
=\int_Z\phi\,d\mu_Z+ O(\min(r^nD)^{-\kappa}\|\phi\|_{C^\theta})
\end{align*}
This shows that 
\begin{align*}
\frac{1}{|D|}\int_{D} f(a m \beta(\alpha^n(\eta(t)))\Lambda)\,dt
&=\int_Z f(az)\,d\mu_Z(z)+O(\min(D)^{-\kappa} r^{-\kappa n}\|f\|_{C^\theta})\\
&=\int_Z f(\alpha^n(h)z)\,d\mu_Z(z)+O(\min(B)^{-\kappa} \rho^{n}\|f\|_{C^\theta})
\end{align*}
with $\rho=r^\kappa\in (0,1)$. This proves (\ref{eq:step}).

Next, we apply the above argument inductively.
For a \holder function $f$ on $X=G/\Lambda$, we define
a function $\bar f$ on  $\bar X=G/M\Lambda$ by 
$$
\bar f(gM\Lambda):=\int_{M/(M\cap\Lambda)} f(gm\Lambda)d\mu_Z(m).
$$
Clearly, 
$$
\|\bar f\|_{C^\theta}\le \|f\|_{C^\theta}.
$$
Let $\bar G=G/M$, $\bar \Lambda=(M\Lambda)/M$,  and $p:G\to \bar G$ be the projection map.
Then $\bar X\simeq \bar G/\bar\Lambda$. We note that $Dp({W}^\alpha)$ is precisely
the unstable space of $D\alpha$ acting on $\mathcal{L}(\bar G)$.
It follows from (\ref{eq:step}) that there exists $\rho\in (0,1)$ such that 
$$
\frac{1}{|B|}\int_B f(\alpha^n(h\psi(b))g\Lambda)\,db=\frac{1}{|B|}
\int_B \bar f(\alpha^n(\bar h\bar \psi(b))\bar g\bar \Lambda)\,db+O(\min(B)^{-\kappa}\rho^n\|f\|_{C^\theta}),
$$
where $\bar\psi$ is the product of the maps of the form
$$
\bar \psi_i: (t_1,\ldots,t_s)\mapsto \exp(t_1\bar w_1)\cdots \exp(t_s\bar w_s),
$$
or 
$$
\bar \psi_i: (t_1,t_1'\ldots,t_s,t_s')\mapsto \exp(t_1\bar w_1+t'_1\bar w'_1)\cdots \exp(t_s\bar
w_s+t_s'\bar w_s').
$$
with $\bar w_i=Dp(w_i)$ and $\bar w'_i=Dp(w_i')$, $\bar h=p(h)$ and $\bar g=p(g)$.
In this product we may skip terms with 
$\bar w_i=0$ or $\bar w_i'=0$ (note that if $\bar w_i=0$, then $\bar w_i'=0$ and conversely).
Then the relations (\ref{eq:jordan1})--(\ref{eq:jordan2}) are still satisfied. 
In particular, the last exponential in the obtained product corresponds to the subspace
$Dp({W}^\alpha)\cap \left<w,\bar w\right>$ where $w$ is an eigenvector of $D\alpha$
in $\mathcal{L}(\bar G)$ with the eigenvalue of maximal modulus.
Now we can again apply the argument as in the proof of (\ref{eq:step}) reducing the number of 
terms in the product defining $\bar \psi$. Repeating the same  argument repeatedly, we deduce that
for some $\rho\in (0,1)$ and $\kappa>0$,
$$
\frac{1}{|B|}\int_B f(\alpha^n(h\psi(b))g\Lambda)\,db=
\int_{M/(M\cap\Lambda)} f(\alpha^n(h)gm\Lambda)\,d\mu_Z(m)+O(\min(B)^{-\kappa}\rho^n\|f\|_{C^\theta}),
$$
where $M$ is a closed normal $\alpha$-invariant
subgroup containing $\exp({W}^\alpha)$ such that $M/(M\cap\Lambda)$ is compact.
We observe that $D\alpha$ acting on $\mathcal{L}(G/M)$ has no eigenvalues with absolute value
greater than one. Since $\alpha$ is ergodic, it follows from Lemma \ref{l:eigenvalue} that $M=G$.
This proves the theorem for the set of $g\in G$ that satisfy (\ref{eq:dense0}) at every inductive step,
with the estimate which is uniform over $g$. Since this set has full measure, we conclude that the 
estimate holds for all $g$ completing the proof of the theorem.
\QED

The following corollary will be used in the proof of the limit theorems in the next section.

\begin{corollary}\label{c:unstable}
Let $\Omega$ be a domain in ${W}^\alpha$ with a piecewise smooth boundary.
Then there exist $\kappa=\kappa(\theta)>0$ and  $\rho=\rho(\theta)\in (0,1)$ 
such that for every $\theta$-H\"older function
$f:X\to \mathbb{R}$, $g\in G$ and $\epsilon>0$, we have
$$
\int_\Omega f(\alpha^n(\psi(b))g\Lambda)\, db=|\Omega|\int_X f\,d\mu+O\left((|\partial_\epsilon\Omega|+\epsilon^{-\kappa}\rho^n|\Omega|)\|f\|_{C^\theta}\right),
$$
where $\partial_\epsilon\Omega$ denotes the $\epsilon$-neighbourhood of the boundary of $\Omega$. 
\end{corollary}

\proof
We tessellate ${W}^\alpha$ by cubes $B$ of size $\epsilon$.
Then
$$
\left|\Omega-\bigcup_{B\subset \Omega} B\right|\le |\partial_\epsilon\Omega|,
$$
and 
$$
\int_\Omega f(\alpha^n(\psi(b))g\Lambda)\, db=\sum_{B\subset\Omega} 
\int_B f(\alpha^n(\psi(b))g\Lambda)\, db+O(|\partial_\epsilon\Omega| \|f\|_{C^0}).
$$
By Theorem \ref{th:unstable}, for some $\kappa>0$ and $\rho\in (0,1)$, 
$$
\int_B f(\alpha^n(\psi(b))g\Lambda)\, db=|B|\int_X f\, d\mu+ O(|B|\epsilon^{-\kappa}\rho^n\|f\|_{C^\theta}).
$$
Therefore,
\begin{align*}
\int_\Omega f(\alpha^n(\psi(b))g\Lambda)\, db &=\left(\sum_{B\subset\Omega} |B|\right)
\int_X f\, d\mu+ O\left(\left(|\partial_\epsilon\Omega|+\sum_{B\subset\Omega}|B|\epsilon^{-\kappa}\rho^n\right)\|f\|_{C^\theta}\right)\\
&=|\Omega|\int_X f\, d\mu+ O((|\partial_\epsilon\Omega|+|\Omega|\epsilon^{-\kappa}\rho^n)\|f\|_{C^\theta}).
\end{align*}
This completes the proof of corollary.
\QED

%%%%%%%%%%%%%%%%%%%%%%%%%%%%%%%%%%
%%%%%%%%%%%%    CLT    %%%%%%%%%%%

\section{Central limit theorem and invariance principles}\label{CLT}

Let us first review the terminology regarding the central limit theorem and other probabilistic limit theorems.
Let $\alpha: X \to X$ be a measure-preserving map of a probability space $(X,\mu)$.
For a function $f: X \to \R$, we consider a sequence of observables $f\circ \alpha^n$.
If the dynamical system $\alpha\curvearrowright X$ is sufficiently chaotic, this sequence is
expected to behave similarly to a sequence of independent random variables.
We set 
$$
S_n (f,x)=\sum _{i=0} ^{n-1} f (\alpha ^i (x)),
$$
and for simplicity assume that $\int_X f\,d\mu=0$.

The sequence $f\circ \alpha^n$ satisfies the {\it central limit theorem}
if there exists $\sigma>0$ such that $n^{-1/2}S_n(f,\cdot)$ converges in distribution
to the normal law with mean $0$ and variance $\sigma^2$.
More generally, the sequence $f\circ \alpha^n$ satisfies the
{\it central limit theorem for subsequences}
if there exists $\sigma >0$ such that for every increasing sequence of measurable functions $k_n (x) $
taking values in $\N$ such that for almost all $x$, $\lim_{n\to\infty} \frac{k_n(x)}{n} =c$ for some
fixed constant $0<c < \infty$,    the sequence  $n^{-1/2} S_{k_n(\cdot)} (f,\cdot)$ converges
  in distribution to the normal law  with mean $0$ and variance $\sigma^2/c$.
We define $S_t(f,x)$ for all $t\ge 0$ by linear interpolation of its values at integral points.
The sequence $f\circ \alpha^n$ satisfies the {\it Donsker invariance principle}
if there exists $\sigma>0$ such that
the sequence of random functions $(n\sigma^2)^{-1/2}S_{nt}(f,\cdot)\in C([0,1])$ converges in distribution
to the standard Brownian motion in $C([0,1])$.
The sequence $f\circ \alpha^n$ satisfies the {\it Strassen invariance principle}
if there exists $\sigma>0$ such that for almost every $x$, the sequence of functions
$(2n\sigma^2\log\log n)^{-1/2}S_{nt}(f,x)$ 
is relatively compact in $C([0,1])$ and its limit set is precisely 
the set of absolutely continuous functions $g$ on $[0, 1]$ such that
$g (0)  = 	0$	and $\int_0^1 g'(t)^2\,  dt \leq 1$.
This is a strong version of the law of the iterated logarithm.
  
In this section we establish the above limit theorems for sequences generated by ergodic
automorphisms of compact nilmanifolds. In the case of toral automorphism,
these theorems have been established by LeBorgne \cite{LB99} using the method
of martingale differences, and we follow a similar approach.
We shall use the following general result:

\begin{Theorem}    \label{Gordin}
Let $(X,{\cal B},\mu,\alpha)$ be an invertible ergodic dynamical system
and $f\in L^2(X)$ such that $\int_X f\, d\mu=0$.
Let ${\cal A}$ be a sub-$\sigma$-algebra of ${\cal B}$ such that 
${\cal A}_n=\alpha^{-n}({\cal A})$ is a non-increasing sequence of $\sigma$-algebras satisfying
\begin{equation}    \label{eq:gordin}
 \sum _{n>0}  \| E(f\mid {\cal A}_n) \| _2  < \infty\quad \text{ and }\quad \sum _{n<0}  \| f - E(f\mid {\cal A}_n) \| _2  < \infty. 
 \end{equation}
Then
\begin{enumerate}
\item[(i)] $\sigma ^2  = \int _X f^2\,  d\mu  +  2  \sum _{j=1}  ^{\infty}   \int _X (  f \circ \alpha ^j )
  f\,  d\mu$ is finite.
\item[(ii)] $\sigma=0$ $\Leftrightarrow$  $f$ is an $L^2$ coboundary $\Leftrightarrow$ $f$ is a measurable coboundary.
\item[(iii)] If $\sigma>0$, then $f\circ\alpha^n$
satisfies the central limit theorem, the central limit theorem of subsequences,
and  the Donsker and Strassen invariance principles.
\end{enumerate}
\end{Theorem}  

It is well-known (see, for instance, \cite[Theorem 4.13]{Viana}) that under the assumption
(\ref{eq:gordin}) the function $f$ has a decomposition $f=(\phi\circ\alpha-\phi)+\psi$ with
$\phi,\psi\in L^2(X)$, where $\psi\circ\alpha^n$ is a reverse martingale difference
with respect to the $\sigma$-algebras $\mathcal{A}_n$, and $\sigma=\|\psi\|_2$.
In particular, $\sigma<\infty$ and if $\sigma=0$, then $f$ is an $L^2$ coboundary.
On the other hand, if $f$ is a measurable coboundary, then $\psi$ is also 
a measurable coboundary, and it follows from \cite{sv} that $\psi=0$, so that $\sigma=0$.
For (iii) we refer to \cite[Ch.~5]{HH}.

The following is the main result of this section:

\begin{Theorem}  \label{affine central limit}
Let $\alpha $ be an  ergodic automorphism of a compact nilmanifold $X$,
and let $f$ be a H\"older function on $X$ which has zero integral
and is not a measurable coboundary.
Then the sequence $f\circ\alpha^n$ satisfies 
the central limit theorem, the central limit theorem of subsequences,
and  the Donsker and Strassen invariance principles.
\end{Theorem}

%Le Borgne and P\`{e}ne  proved a multiple central limit theorem in \cite{LeBorgne-Pene}.  We expect that our methods will permit a generalization of her results to nilautomorphisms  but do not pursue this here.

To  find the sub-$\sigma$-algebra $\mathcal{A}$ suitable for Theorem \ref{Gordin}, 
we use the results of Section \ref{sec-unstable} combined with the works of Lind \cite{Lind} and  Le
Borgne \cite{LB99}.  
We call a measurable partition $\mathcal{P}$ of $X$ {\it $\delta$-fine} if the diameter of any set in $\mathcal{P}$ is
at most $\delta$. We say that a partition generates under $\alpha$ if the $\sigma$-algebra generated by all
$\alpha ^n (\mathcal{P})$ with $n \in \Z$ is the Borel $\sigma$-algebra of $X$ modulo null sets. 
Given a partition $\mathcal{P}$ and $x \in X$, we denote by $\mathcal{P} (x)$ the element of the
partition that contains $x$.  Given  integers $k \leq l$, we denote by $\mathcal{P}^l _k$ the partition
generated by $\alpha ^{-k} (\mathcal{P}), \ldots, \alpha ^{-l} (\mathcal{P})$.
We also set $\mathcal{P}^{\infty} _k (x) = \cap _{l \geq k} \mathcal{P}^l _k (x)$.

%The following result for toral automorphisms was proved by Lind \cite{Lind}.  

%Next we construct a suitable partition $\mathcal{P}$ which allows to construct 
%$\sigma$-algebras $\mathcal{A}_n$ as in Theorem \ref{Gordin}.

\begin{Proposition}\label{p:partition}
Let   $\mathcal{P}$ be a finite measurable partition of $X$ such that for every $P \in \mathcal{P}$,
\begin{itemize}
\item $P$ is the closure of its interior,  
\item the boundary of $P$ is piecewise smooth,
\item the diameter of $P$ is at most $\delta$.  
\end{itemize}
Then if $\delta$ is sufficiently small,
\begin{enumerate}
\item[(i)] the partition $\mathcal{P}$ generates under $\alpha$,
\item[(ii)] for almost every $x$, the atoms  $\mathcal{P}^{\infty} _0 (x)$ are contained in the stable manifolds
$\mathcal{W}^{s} (x)$ of $x$,
and the diameter of  $\mathcal{P}^{\infty} _0 (x) $  in $\mathcal{W}^{s} (x)$ is  bounded,
\item[(iii)] for almost every $x\in X$, the atoms $\mathcal{P}^{\infty} _0 (x) $ have non-empty interior  in the stable manifolds
$\mathcal{W} ^{s} (x)$. 
\end{enumerate}
\end{Proposition}

 {\em Proof of (i)--(ii).}  The proof follows that of \cite[Th.~1]{Lind} almost completely albeit with some differences
in the final argument involving isometries.   We will show that $\alpha$ almost surely separates points, i.e., that for
some null set $X_0$ in $X$, if $x,y  \in  X\backslash X_0$, then for some $n$, the points $\alpha ^n (x)$ and $\alpha ^n (y)$
belong to different elements of the partition $\mathcal{P}$.  It then follows from Rohklin's work \cite{Rok} that
$\mathcal{P}$  generates under $\alpha$.  

There exist $c_0>1$ and $\delta_0>0$ such that for every $w\in\mathcal{L}(G)$ 
satisfying $\|w\|<\delta_0$ and $x\in X$,
\begin{equation}\label{eq:compare1}
c_0^{-1}\,\|w\|\le d(x,\exp(w)x)\le c_0\, \|w\|.
\end{equation}
We assume that $\delta$ is sufficiently small, so that 
$\|D\alpha\|c_0\delta<\delta_0$, and
if $p$ and $q$ belong to the same
element $P$ of the partition, then $q=\exp(w)p$ with $\|w\|<\delta_0$.
Since $\hbox{diam}(P)\le \delta$, we have $\|w\|\le c_0\delta$. 
We observe that
\begin{equation}\label{eq:compare2}
d(\alpha^n(p),\alpha^n(q))=d(\alpha^n(p), \exp((D\alpha)^n w)\alpha^n(p)).
\end{equation}
Suppose that $\|(D\alpha)^n w\|\to \infty$ as $n\to\infty$.
We pick the greatest  $n\ge 0$ such that $\|(D\alpha)^n w\|\le c_0\delta$.
Then
$$
c_0\delta< \|(D\alpha)^{n+1} w\|\le \|D\alpha\|c_0\delta<\delta_0,
$$
and it follows from \eqref{eq:compare1}--\eqref{eq:compare2} that
$d(\alpha^{n+1}(p),\alpha^{n+1}(q))>\delta$. Hence, 
$\alpha^{n+1}(p)$ and $\alpha^{n+1}(q)$ belong to different elements of the partition.

A similar argument also applies when $\|(D\alpha)^n w\|\to \infty$ as $n\to -\infty$.
Therefore, it remains to consider the case when $w\in E^{iso}$ which is the span of eigenspaces of 
$D\alpha$ with eigenvalues of modulus one.
We adapt Lind's idea \cite{Lind} for this
situation.  Let $K$ denote the closed group of isometries generated by $\beta :=D \alpha
|_{E^{iso}}$.  
Then $\beta$ acts ergodically on $K$ by translations. Since $\alpha$ is mixing, the product $\alpha \times
\beta$ acts ergodically on $X \times K$. It follows from ergodicity
and  Fubini's theorem that there exists a null set $X_0\subset X$ and $k\in K$ such that
the sequence $(\alpha^n(x), \beta^nk)$ is dense in $X\times K$
for every $x\in X\backslash X_0$.
Then the sequence $(\alpha^n(x), \beta^n)$ is also dense in $X\times K$.

Now suppose that $p,q\in X\backslash X_0$ and $q=\exp(w)p$ for some nonzero $w\in E^{iso}$.
Given an element $P\in \mathcal{P}$, we set
$$
P(w,\epsilon)=\{x\in P:\, d(\exp(w)x,P)>\epsilon \}.
$$
When $\epsilon>0$ is sufficiently small, this set has a nonempty interior.
Hence, for every $p\in X\backslash X_0$, there exists $n$ such that
$$
\alpha^n(p)\in P(w,\epsilon)\quad\hbox{and}\quad d(\exp(w),\exp((D\alpha)^n w))<\epsilon/2.
$$
Then 
\begin{align*}
d(\alpha^n(q),P)&=d(\exp((D\alpha)^n w) x,P)\\
&\ge
d(\exp(w) x,P)-d(\exp((D\alpha)^n w)x,\exp(w)x)>\epsilon/2.
\end{align*}
In particular, $\alpha^n(p)\in P$ and $\alpha^n(q)\notin P$. 
This proves that $\mathcal{P}$ generates under $\alpha$.
The part (ii) can be proved by the same argument.
\QED

To prove Proposition \ref{p:partition}(iii), 
we follow Le Borgne's approach [17] for toral automorphisms.
We pick $c, r_0 \in (0,1)$ such that the map $\alpha ^{-n}$ expands 
the distance on $\mathcal{W} ^{s}$ by at least $c\, r_0 ^{-n}$ for $n \geq 0$,
and take $r\in (r_0,1)$. Let
\[ V_n  := \{x\in X:\, \mathcal{P} _0 ^{\infty} (x )\supset
B_{  r ^n/ c }  (x)\cap \mathcal{W}^s(x) \}  .\]
Proposition \ref{p:partition}(iii) immediately follows from the following lemma.

\begin{Lemma}  \label{lem-boundary}
$\mu (X \backslash V _n )  \ll  r ^n$. 
\end{Lemma} 

\proof
Let  
$$
W_n:=\{ y\in X:\, d(\alpha ^j (y), \partial \mathcal{P}(\alpha ^j (y)))  \geq  r_0 ^j  r ^n / c ^2 \text{
  for all } j \geq 0\}.
$$ 
 If $y$ is in $W_n$, then $\mathcal{P} (\alpha ^j (y))$ contains the ball in $\mathcal{W} ^s (\alpha ^j
 (y))$  of radius $r_0 ^j  r ^n /c ^2$.  Hence, $\alpha ^{-j} (\mathcal{P}(\alpha ^j (y)))$ contains the
 ball in $\mathcal{W} ^s (y)$ of radius $r ^n /c$. Since
$$
\mathcal{P} _0 ^{\infty} (y) = \bigcap _{j \geq 0} \alpha ^{-j}(\mathcal{P}(\alpha ^j (y))),
$$
we conclude that $V_n \supset W_n$.  
 
 To prove the lemma, it suffices to estimate $\mu (X\backslash W_n)$.
It follows from our assumption on the partition $\mathcal{P}$ that
$$
\mu ( \{y\in X:\,  d(y, \partial
 \mathcal{P}(y))  \leq  \epsilon \})  \ll \epsilon,
$$
and since $\alpha$ is measure-preserving, for every $j\ge 0$,
$$
\mu ( \{y\in X:\,  d(\alpha ^j (y), \partial
 \mathcal{P}(\alpha ^j (y)))  \leq  r_0 ^j  r ^n  / c\})  \ll r_0 ^j r ^n.
$$
Hence,
$$
\mu(X\backslash W_n) \ll \sum _{j\geq 0} r_0 ^j r ^n \ll r ^n,
$$   
which implies the lemma.
\QED

%\vspace{1em}

{\em Proof of Theorem \ref{affine central limit}.}
Let $\mathcal{A}$ be the $\sigma$-algebra generated by the partition  $\mathcal{P}^{\infty} _0$ and
$\mathcal{A}_n = \alpha^{-n} (\mathcal{A})=\mathcal{P}^{\infty} _n$.
It is clear that the sequence $\mathcal{A}_n$ is non-increasing.
To prove the theorem, it suffices to check the conditions (\ref{eq:gordin}).
Since the partition $\mathcal{P}^{\infty} _n$  is measurable  in the sense of \cite{Rok},
for almost every $x$,
\[
E(f|\mathcal{A}_n) (x)  = \int _{\mathcal{P}^{\infty} _n (x)} f(y)\,  dm_{\mathcal{P}^{\infty} _n (x)}(y),
\]
where $m_{\mathcal{P}^{\infty} _n (x)}$ is the conditional probability measure on $\mathcal{P}^{\infty}
_n (x)$. 

To verify the second part of (\ref{eq:gordin}), we observe that when 
$\mathcal{P}_0^\infty(\alpha^n(x))\subset \mathcal{W}^s(x)$,
$$
\hbox{diam}(\mathcal{P}_n^\infty(x))=\hbox{diam}(\alpha^{-n}(\mathcal{P}_0^\infty(\alpha^n(x))))
$$
decays exponentially as $n\to -\infty$ uniformly on $x$. Since the function $f$ is $\theta$-H\"older,
it follows that for some $\tau\in (0,1)$,
$$
\| f - E(f\mid {\cal A}_n) \| _2\ll \tau^{-n}\|f\|_{C^\theta}\quad\hbox{and}\quad
\sum _{n<0}  \| f - E(f\mid {\cal A}_n) \| _2  < \infty.
$$

To check the other condition in (\ref{eq:gordin}), we observe that by
Lemma \ref{lem-boundary}, 
\begin{equation}\label{eq:complement}
   \int _{X \backslash \alpha^{-n}(V_n)} |E(f|\mathcal{A}_n)|^2\, d \mu  \ll  r ^n  \| f \|^2 _{C^{0}}.
\end{equation}
  On the other hand, for  $x \in \alpha^{-n}(V_n)$,
 \[    B_{  r ^n/ c }  (\alpha ^{n} (x))\cap \mathcal{W}^s(\alpha^n(x))   \subset  \mathcal{P} (\alpha
 ^{n} (x))\quad\hbox{and}\quad
B_{  r ^n r_0 ^{-n} } (x)\cap\mathcal{W}^s(x) \subset \alpha ^ {-n}(\mathcal{P} (\alpha ^{n} (x))).\]
 Since the diameter of $\mathcal{P} (x)$ is at most $\delta$, as soon as $ r^n r_0^{-n} > \delta$, we get that $\mathcal{P} (x) \subset B_{  r^n/r_0^n } (x) $.
 Hence, by Proposition \ref{p:partition}, for almost every $x\in\alpha^{-n}(V_n)$,
 \begin{equation}\label{eq:finite intersection}
\mathcal{P} _0 ^{\infty} (x) = 
\bigcap _{n=0} ^{\infty} \alpha ^{-n} (\mathcal{P} (\alpha ^{n} (x))) \cap\mathcal{W}^s(x)=\bigcap _{n=0} ^{\left\lceil\frac{\log \delta} {  \log (r/ r_0)}\right\rceil} \alpha ^{-n} (\mathcal{P} (\alpha ^{n} (x))) \cap\mathcal{W}^s(x).
 \end{equation}
 Thus, $\mathcal{P}_0^{\infty} (x)$ is the intersection of the stable manifold of $x$  with at most
finitely many sets  whose boundaries consist of finitely many piecewise smooth submanifolds.
Then
\begin{equation}\label{eq:p_0}
\mathcal{P} _0 ^{\infty} (x) =\exp(\Omega_x)x,
\end{equation}
where $\Omega_x$ is a domain in the unstable subspace $W=W^{\alpha^{-1}}$ of $D(\alpha^{-1})$ in $\mathcal{L}(G)$
whose boundary is piecewise smooth and depends smoothly on $x$. In particular,
$|\partial_\epsilon \Omega_x|\ll \epsilon$
uniformly on $x\in X$. It follows from (\ref{eq:p_0}) that
$$
\mathcal{P} _n ^{\infty} (x)=\alpha^{-n}(\mathcal{P} _{0} ^{\infty} (\alpha^n(x))) =\exp((D\alpha)^{-n}\Omega_x)x.
$$
Then by \cite[Prop.~4.3]{CB0}, 
$$
m_{\mathcal{P}^{\infty} _n (x)}=\frac{1}{m_x(\mathcal{P}^{\infty} _n (x))}m_x|_{\mathcal{P}^{\infty} _n
  (x)},
$$
where $m_x$ is the Haar measure on $\exp(W)x$.
Now we apply Corollary \ref{c:unstable}.  It follows from the definition of $V_n$ that for $x\in \alpha^{-n}(V_n)$,
we have $|\Omega_x|\gg r^{n}$. Hence, by Corollary \ref{c:unstable},
for every $x\in \alpha^{-n}(V_n)$ and $\epsilon>0$,
\begin{align*}
\frac{1}{m_x(\mathcal{P}^{\infty} _n (x))} \int_{\mathcal{P}^{\infty} _n (x)}f(y)\, dm_x(y)
&=O\left(\left(\frac{|\partial_\epsilon\Omega|}{|\Omega|}+\epsilon^{-\kappa}\rho^n\right)\|f\|_{C^\theta}\right)\\
&=O\left((\epsilon r^{-n}+\epsilon^{-\kappa}\rho^n)\|f\|_{C^\theta}\right),
\end{align*}
where $\rho\in (0,1)$. We take $\epsilon=(r^n\rho^n)^{1/(\kappa+1)}$.
If we also take $r$ sufficiently close to $1$, then this quantity decays exponentially as $n\to\infty$.
Then
$$
\int _{\alpha^{-n}(V_n)} |E(f|\mathcal{A}_n)|^2\, d \mu  \ll  \tau ^n  \| f \|^2 _{C^{\theta}}
$$
for some $\tau\in (0,1)$. Combining this estimate with (\ref{eq:complement}),
we deduce the first part of (\ref{eq:gordin}). Now the theorem follows from Theorem \ref{Gordin}.
\QED

% Given a subset $P \subset {\cal W}^s (x)$, let the $size (P)$ be the radius of the largest inscribed ball (relative to a fixed ambient Riemannian metric).  Furthermore, denote by $\partial _{\varepsilon} $  the $\varepsilon$-neighborhood of a set in ${\\|         cal W}^s$.   We observe that 

%\begin{enumerate}
%\item $\mu \{x \mid  size ({\cal A} (x) < \varepsilon )  = O(\varepsilon)$
%\item $ \mu ^s ( \partial _{\varepsilon} {\cal A} (x) )  =  O(\varepsilon)$.
%\end{enumerate}

%Let us now verify the two  conditions in  Theorem \ref{Gordin}.   The second condition follows since the ${\cal A}_n = \alpha ^{-n} ({\cal A})$ are shrinking uniformly as $n \rightarrow -\infty$.  Hence  $E(f\mid {\cal A}_n)$ approach $f$ with uniform estimates.   The       first condition is the hard one, and uses equidistribution of unstable manifolds.  We have the following 

%%%%%%%%%%%%%%%%%%%%%%%%%%%%%%%%%%
%%%%%%%%%%%%    COHOMOLOGY    %%%%%%%%%%%

\section{Cohomological equation}\label{sec:cohomology}

In this section we apply exponential mixing to establish regularity of solutions of 
the cohomological equation. We recall that 
for ergodic systems the solution is unique up to a constant, up to measure zero.
 
\begin{Theorem} \label{cohomology}
Let $\alpha$ be an ergodic automorphism of a compact nilmanifold $X$  
and $f\in C^\infty(X)$ such that $f=\phi\circ \alpha-\phi$ for some measurable function $\phi$.
Then $\phi$ is almost everywhere  equal to a $C^\infty$ function.
\end{Theorem}

The method of  proof of Theorem \ref{cohomology}
applies to other classes of homogeneous partially hyperbolic systems for which exponential mixing 
 holds.
For instance, we may consider an ergodic partially hyperbolic left translation on
the homogeneous space $G/\Gamma$, where $G$ is connected semisimple Lie group and $\Gamma$
is a cocompact irreducible lattice. This dynamical system is also 
exponentially mixing for H\"older functions \cite[Appendix]{Kleinbock-Margulis}, and the argument
of Theorem \ref{cohomology} applies. For $X=\hbox{SL}_d(\mathbb{R})/\hbox{SL}_d(\mathbb{Z})$,
an analogous result for H\"older functions $f$ was established in \cite{LB02}. Furthermore, we get both \holder and smooth versions of Theorem \ref{cohomology} 
for compact $G/\Gamma$ and $G$ semi simple from  Wilkinson's general result for accessible partially hyperbolic diffeomorphisms \cite[Theorem A]{Wilkinson}  under the additional assumption  that the left translation projected to any factor of $G$ does not belong to a compact subgroup.

Before starting the proof, we need to develop some language and review a result on regularity of
distributions. Let $M$ be a compact manifold.  
We fix a Riemannian metric on $M$, and denote by  $C^{\theta}=C^\theta(M)$  the space of   $\theta $-H\"{o}lder
 functions on $M$.  We let $(C^{\theta})^*$ be the dual space to $C^{\theta} $.  Note that any   smooth
 function on $M$ naturally belongs to any $C^{\theta}$.  Hence any element in $(C^{\theta})^*$ defines a
 distribution on smooth functions on $M$.  Conversely, $(C^{\theta} )^*$ is the space of distributions
 (dual to \smooth functions) which extend to continuous
 linear functionals on  $C^{\theta}$.
 As for notation, we will  write the pairing  $D(g) = \langle D, g \rangle$ for $D \in (C^{\theta} )^*$ and $g \in C^{\theta} $.  
 
 Let  ${\cal F}$ be a \smooth  foliation  on $M$,  and consider a \smooth vector field $V$ tangent to ${\cal F}$.  Given a distribution $D$ on $M$,  define the derivative $V(D)$ by evaluating on \smooth test functions $g$ as follows:  $\langle V(D), g \rangle = - \langle D, V(g) \rangle$ where $V(g)$ denotes the directional derivative of $g$ along $V$.    
 % Thus for $D$ to have partial derivatives of order $m$ to belong to $(C^{\theta} )^*$ is well defined and independent of the choice of spanning vector fields.  

 Given  smooth vector fields $V_1, \ldots, V_r$, we call $V_{i_1}, V_{i_2} \ldots V_{i_m} D$ the partial derivatives of order $m$ of $D$.  Suppose that we can cover $M$ with open sets $\mathcal{U}$ such that we 
can find smooth vector fields $V_1, \ldots, V_r$ which span the tangent spaces to ${\cal F}$ at any point of $\mathcal{U}$.  Suppose moreover that all partial  derivatives of any order $m$, $V_{i_1}, V_{i_2} \ldots V_{i_m} D$ of a distribution $D$  belong to  $(C^{\theta} )^*$, for all such choices of $\mathcal{U}$ and $V_1, \ldots V_r$ .  Then  for any other \smooth vector fields $V'_1, \ldots V'_r$ tangent to ${\cal F}$, the partial derivatives 
$V'_{i_1}, V'_{i_2} \ldots V'_{i_m} D$ also belong to   $(C^{\theta} )^*$ as follows from a partition of
unity argument.   Thus we can say that partials along ${\cal F}$ of a distribution  belong to  $(C^{\theta} )^*$, without any reference to a particular set of  vector fields.\footnote{In our application we will have globally defined vector fields for which the partials exist for all orders, and we will not need this  comment.}

The following result is inspired by results of Rauch and Taylor in \cite{RT}, and was known to Rauch for the case of \smooth foliations. We are not aware of a simple reference.  It is also   a straight-forward consequence of a similar much more technical result for \holder foliations proved in \cite{FKS}, namely  that the wavefront set  of a distribution for which the partial derivatives of all orders along a single foliation belong to the dual of \holder functions  is co-normal to the foliation. We refer to \cite{RT,FKS} for more details.

\begin{corollary}[\cite{FKS}] \label{what we need}
Let  ${\cal F}_1, \ldots, {\cal F}_r$ be \smooth  foliations  on a  compact manifold $M$ whose  tangent spaces span the tangent spaces to $M$ at all points.
  Consider  a distribution $D$ defined by integration against  an $L^1$ function $\phi$.   Assume that any
  partial derivative of $D$ of any order along the foliations  ${\cal F}_1, \ldots, {\cal F}_r$ belongs to $(C^{\theta})^*$ for all $\theta>0$.
Then $\phi$ is $C^\infty$. 
\end{corollary}

We are now ready to tackle the proof of Theorem \ref{cohomology}.
 Let us first give an outline of the argument.  Using Theorem \ref{Gordin}, we  first show in Lemma
 \ref{lemma-holder} that the function $\phi$ has to be in $L^2(X)$.  Then we describe $\phi$ as
 distribution.
We consider three dynamically defined  foliations for $\alpha$:
the unstable foliation $\mathcal{W}^u$, the stable foliation $\mathcal{W}^c$, and 
the central foliation $\mathcal{W}^c$. The unstable foliation 
is tangent to the
right invariant distribution on $X$ corresponding to the sum of all generalized eigenspaces with eigenvalues
$|\lambda|>1$, the stable foliation is tangent to the
right invariant distribution on $X$ corresponding to the sum of all generalized eigenspaces with eigenvalues 
$|\lambda|<1$, and the central foliation 
is tangent to the
right invariant distribution on $X$ corresponding to the sum of all generalized eigenspaces with eigenvalues 
$|\lambda|=1$.
Note that these  distributions are integrable as is easily seen by taking Lie brackets.
We show that the distribution derivatives of $\phi$ along the foliations 
$\mathcal{W}^s$, $\mathcal{W}^u$, $\mathcal{W}^c$ of
 $\alpha$ define distributions on \holder functions. This is established in Lemmas
 \ref{lemma:derivativesofhstable} and \ref{lemma:derivativesalongneutral}.  Since  all these foliations
 are smooth, 
  Corollary \ref{what we need} shows that the function $\phi$ is $C^{\infty}$. 

We now establish Lemmas \ref{lemma-holder}, \ref{lemma:derivativesofhstable} and \ref{lemma:derivativesalongneutral} which will finish the proof of Theorem \ref{cohomology}. 
\begin{Lemma}  \label{lemma-holder}  
The function $\phi$ in Theorem \ref{cohomology} is in $L^2$.
\end{Lemma} 

\proof Recall that along the proof of Theorem \ref{affine central limit} we have verified the conditions
of Theorem \ref{Gordin}. Hence, the lemma follows from part (ii) of this theorem.
 \QED
  
 Define the distributions $P^+$ and $P^-$ by evaluating them on test functions $g\in C^\infty(X)$ by
  \[  P^+ (g)   = \sum _{i=0} ^{\infty} \langle f \circ \alpha ^i , g
 \rangle   \text{ \hspace{2em} and \hspace{2
 em}}  P^- (g)   = \sum _{i=1} ^{\infty} \langle f \circ \alpha ^{-i} , g \rangle.  \]
Note that $\int _X f\, d\mu =0$ since $f$ is an $L^2$ coboundary.
Hence, by exponential mixing (Theorem \ref{exp mixing}), these sums converge as long as the test function $g$ is
  H\"{o}lder,
and $P^+,P^-\in (C^\theta)^*$.  
Moreover, since $ \langle \phi \circ \alpha ^i , g \rangle  \rightarrow 0$ as $i \rightarrow \pm \infty$.   we get by a telescoping-sum argument that 
$$P^+ (g) =  \sum _{i=0} ^{\infty} \langle f \circ \alpha ^i , g \rangle   =  \sum _{i=0} ^{\infty} \langle
\phi \circ \alpha ^{i+1} - \phi \circ \alpha ^i  , g \rangle = \lim _{N \rightarrow \infty} \langle \phi \circ \alpha ^N -\phi , g \rangle  = -
\langle \phi , g \rangle.$$
Similarly, we see that $P^-(g)=\langle \phi , g \rangle$. Hence, the distribution 
$P^{+} = -  P^{-}$  is given by integration against the $L^2$-function $\phi$.  We will use this to show that $\phi$ is smooth.  

According to Corollary \ref{what we need}, it suffices to show that partial derivatives of all orders of the distribution $P^+=-P^-$   along any of the three foliations ${\cal W}^s$,   ${\cal W}^u$ and ${\cal W}^c$ belong to $(C^{\theta})^*$ for any $\theta>0$.  We will show this in the next two lemmas.
%First,  partials of $P^+$ along ${\cal W}^s$ are always dual to \holder functions due to the dynamics.    The fact that $P^+ = - P^-$ proves the same for partial derivatives along ${\cal W}^u$.  Note that we  use  exponential mixing to insure  that $P^+= - P^-$ are defined on \holder functions. 
 
 \begin{lemma}
\label{lemma:derivativesofhstable}
 Partial derivatives of all orders of the distribution $P^+=-P^-$   along  ${\cal W}^s$ and ${\cal W} ^u$ belong to $(C^{\theta})^*$ for any $\theta>0$.
\end{lemma}

\proof Let $V$ be  a right  invariant vector field tangent to $\mathcal{W}^s$ and $g$ a \smooth test function.  Then
\[  \langle V(P^+), g\rangle  =  -  \langle P^+, V(g)\rangle
= - \sum _{i=0} ^{\infty} \langle f \circ \alpha ^i , V(g)
 \rangle
=\sum _{i=0} ^{\infty} \langle V(f \circ \alpha ^i) , g
 \rangle.
\]
The derivative $V(f \circ \alpha ^i)$ decays exponentially fast since $V$ is tangent
to $W^s$. Hence,
\[  |\langle V(P^+), g\rangle|\ll \|g\|_{C^0}, \]
and in particular, $V(P^+)\in (C^\theta)^*$ for all $\theta>0$.
Since $P^+=-P^-$, an analogous proof shows that $V(P^+)$  lies in the dual of \holder functions for all vector fields $V$ tangent to $\mathcal{W}^u$.  
A similar argument also applies to higher order derivatives along vector fields tangent to
$\mathcal{W}^s$ or $\mathcal{W}^u$. 
We refer for the details  to \cite[Lemma 5.1]{FKS}.
\QED

Finally, we  show that partials of all orders of $P^+=-P^-$ along ${\cal W}^c$ are distributions on \holder functions.  This argument uses exponential decay very strongly, and was first discovered  in \cite{FKS}.  For a detailed account we refer to \cite[Lemma 5.1]{FKS}.

\begin{lemma}
\label{lemma:derivativesalongneutral}
 Partial derivatives of all orders of the distribution $P^+=-P^-$   along  ${\cal W}^c$ belong to $(C^{\theta})^*$ for any $\theta>0$.
\end{lemma}

\proof  Let $V$ be a right invariant vector field  tangent to $\mathcal{W}^c$, and let $g$ be a \smooth function.
Then the  partial derivative of $P^+$ along  $V$ is given by 
\begin{equation}    \label{solution}
\langle V( P^+), g \rangle=  \sum _{i=0} ^{\infty} \langle V ( f \circ \alpha ^i ) , g \rangle 
= - \sum _{i=0} ^{\infty}  \langle f \circ \alpha ^i  , V (g) \rangle,
\end{equation}
and we  have estimates for
all of these expressions in terms of the \holder norm of $V(g)$, due the exponential mixing of $\alpha$.
We will show that  this distribution extends to \holder functions $g$ by approximating  $g$ by smooth
functions $g_{\varepsilon}$ and  carefully balancing the speed of the approximation with the loss of
exponential decay due to the growth of the $C^l$-norm of $g_{\varepsilon}$.
More precisely, we shall show that there exists $\xi=\xi(\theta)\in (0,1)$
such that for every $g\in C^\theta(X)$ and sufficiently large $i$,
\begin{equation}\label{eq:hhh}
| \langle V(f \circ \alpha ^i) , g \rangle |
\ll \xi^{i} \cdot \|f\|_{C^1} \|g\|_{C^\theta}.
\end{equation}
It would follow from (\ref{solution}) and (\ref{eq:hhh}) that $V(P^+)\in (C^\theta)^*$.

We recall from  Lemma \ref{convolution-approximation} that for $\varepsilon >0$,  there is a \smooth function $g_{\varepsilon}$ such that 
\begin{equation} \label{fe}
 \|g_{\varepsilon } -g \|_{C^0} \leq \varepsilon ^{\theta} \|g\|_{C^\theta} \;\;\;
   \text{and} \;\;\;
\|g_{\varepsilon}\|_{C^2}  \ll \, \varepsilon ^{-m-2}  \|g\| _{C^0}
\end{equation}
where $m=\dim(X)$. We first estimate $| \langle V(f \circ \alpha ^i) , g_\epsilon \rangle |$.
  By the exponential mixing (Theorem \ref{exp
  mixing}) and since $V$ is bounded,  we have for some $\rho\in (0,1)$,
\begin{equation}  \label{pfa}
| \langle f \circ \alpha ^i  , V (g_\e) \rangle  |   \ll
     \rho ^i \|f\|_{C^1}   \|V(g_\e) \| _{C^ 1}    \ll \rho ^i \|f\|_{C^1} \, \varepsilon ^{-m-2}  \|g\| _{C^0}.
\end{equation}
On the other hand, we can estimate $|\langle f \circ \alpha ^i , V(g-g_\e) \rangle|=|\langle V(f \circ \alpha ^i) , g-g_\e \rangle|$ 
as follows.  First, we note that the derivatives  $V(f \circ \alpha ^i)$ by the chain rule are composites
of derivatives of $f$ and derivatives of $\alpha ^i$ along  ${\cal W}^c$.  The latter grows at most
polynomially  since  ${\cal W}^c$ is the central foliation for $\alpha$. Hence,  for any $\eta>0$, there
is $i _{\eta}\in \Z _+$ such that 
$$
\| D(\alpha^i|_{\mathcal{W}^c}) \| < (1+\eta)^{i}\quad\hbox{ for all $i \geq i _\eta$.}
$$
Hence, for all  $i > i _\eta$, we get the estimate 
$$
 \| V(f \circ \alpha ^{i}) \| _{C^0} \leq (1+\eta)^{ i} \|f\|_{C^1},
$$
and
\begin{equation}  \label{pf-fa}
|    \langle V(f \circ \alpha ^{i}) , g - g_\e \rangle |\le 
 \| V(f \circ \alpha ^{i}) \| _{C^0} \: \| g - g_\e \| _{C^0}    \le  (1+\eta)^{i} \|f\|_{C^1}\,  \varepsilon ^{\theta} \|g\|_{C^\theta}.
\end{equation}
%Hence we can assume that for any $\eta>0$, $ \| V(f \circ \alpha ^i) \| _0 \leq (1+\eta)^{i}$ for all large $i$, and we get from the last equation that for some $C>0$ 
%\begin{equation} 
%|   \langle V( P^+), g - g_\e \rangle  |   \le \sum _{i=0} ^{\infty}  \| V(f \circ \alpha ^i) \| _0 \: \varepsilon ^{\theta} \|g\|_{\theta} \quad  \le \sum _{i=0} ^{\infty}  C (1+\eta)^{i}\: \varepsilon ^{\theta} \|g\|_{\theta} \quad 
%\end{equation}
We have  exponential
decay with respect to $i$ in \eqref{pfa}, but    exponential growth in \eqref{pf-fa}  at first sight.
However, choosing $\e$ carefully depending on $i$, we can still achieve exponential decay in 
 \eqref{pf-fa}, and hence for $|\langle V ( f \circ \alpha ^{i }) , g \rangle|$. More precisely,
we take
 $
\varepsilon = \rho ^{\frac{i }{\theta + m +2}}
$.
Then we obtain from \eqref{pfa} that
$$
| \langle f \circ \alpha ^i  , V (g_\e) \rangle  |   \ll
 \left(\rho ^{\frac{\theta}{\theta + m +2} } \right) ^i  \:\: \|f\|_{C^1} \,\|g\| _{C^0},
$$
and from  \eqref{pf-fa} that for $ i> i_\eta$,
$$
 |    \langle V(f \circ \alpha ^{i}) , g - g_\e \rangle |\le 
  \left((1+ \eta)  \rho^{\frac{\theta}{\theta + m +2}} \right)^i     \|f\|_{C^1} \|g\|_{C^\theta}.
$$
Now we choose $\eta>0$ so that
$    \xi :=  (1+ \eta)  \rho^{\frac{\theta}{\theta + m +2}} <1       $.    Finally, we obtain from the last two inequalities that for $ i >i_\eta$,
$$
| \langle V(f \circ \alpha ^i) , g \rangle |
\ll \xi^{i} \cdot \|f\|_{C^1} \|g\|_{C^\theta}.
$$
This proves (\ref{eq:hhh}) and shows that
$V(P^+)$ extends to a continuous linear  functional on the space
of $\theta$-H\"{o}lder functions.  
A similar argument shows that higher order derivatives of $P^+$ along the central foliation define distributions dual to \holder functions. For the details we refer to \cite[Lemma 5.1]{FKS}.     \QED

%\vspace{0.5cm}

This finishes the proof of Theorem \ref{cohomology}.

%%%%%%%%%%%%%%%%%%%%%%%%%%%%%%%%%%
%%%%%%%%%%%%    BERNOULLI    %%%%%%%%%%%

\section{Bernoulli property}\label{s:ber}

Here we show that ergodic automorphisms on compact nilmanifolds are Bernoulli combining results from
\cite{Katz-bernoulli}, \cite{Marcuard}, and \cite{Rudolph78}.
It was already shown in \cite{Parry} that such automorphisms satisfy the Kolmogorov property.

\begin{Theorem} \label{Bernoulli}
Ergodic automorphisms on compact nilmanifolds are Bernoulli.
\end{Theorem}

\Proof 
% Recall the result of Ornstein that roots of Bernoulli automorphisms are Bernoulli \cite[Theorem 4.29]{Walters}.  Hence it suffices to show that ergodic automorphisms of infra-nilmanifolds are Bernoulli.  Further, such an automorphism is finitely covered by an automorphism of a nilmanifold.  By  \cite[Theorem 4.29]{Walters}, 
Let $\alpha$ be an ergodic automorphism of a compact nilmanifold $X=G/\Lambda$.
We will argue by induction on the dimension of $X$.
We note that when $X$ is a torus, this result was established by Katznelson in \cite{Katz-bernoulli},
and this forms the base of induction.
Let $Z$ be the centre of $G$. It follows from \cite[5.2.3]{CG} that $Z\Lambda$ is a closed subgroup of $G$.
Then $\alpha\curvearrowright X$ is measurably isomorphic to a skew product with the base $\alpha\curvearrowright Y=G/(Z\Lambda)$ and fibers isomorphic to 
the torus $T=Z\Lambda/\Lambda$, where the action on the fibers is by affine linear maps
$t\mapsto z_y+\alpha(t)$, $z_y\in Z$. We consider two cases.

First, suppose that the automorphism $\alpha$ acts ergodically on the torus $T$.
Then it follows from Marcuard's theorem  \cite[Theorem 4]{Marcuard} 
that $\alpha\curvearrowright X$ 
is measurably isomorphic to the direct product of the systems $\alpha\curvearrowright Y$
and $\alpha\curvearrowright T$.  Hence,  
it follows from the inductive assumption that $\alpha\curvearrowright X$ is measurably isomorphic to the product
of two Bernoulli maps, and thus Bernoulli. % (see \cite[Examples after Definition 4.11]{Walters}). 

Second, suppose that the action of $\alpha$ on the torus $T=Z\Lambda/\Lambda$ is not ergodic. 
Then $T$ contains a nontrivial subtorus $T_0=Z_0\Lambda/\Lambda$,
where $Z_0$ is a closed connected subgroup of $Z$, on which $\alpha$ acts isometrically,
and $\alpha\curvearrowright X$ is measurably isomorphic to a skew product with the base
$\alpha\curvearrowright G/(Z_0\Lambda)$ and the fibers
isomorphic to torus $T_0$, where the action on the fibers is by affine linear maps
$t\mapsto z_y+\alpha(t)$, $z_y\in Z_0$.
We note that this is an isometric extension of the base, and by the inductive assumption, the
base is Bernoulli.   Hence, we can apply Rudolph's theorem \cite{Rudolph78} which shows that weakly mixing isometric
extensions of Bernoulli maps are Bernoulli.  \QED

%%%%%%%%%%%%%%%%%%%%%%%%%%%%%%%%%%
%%%%%%%%%%%%    BIBLIOGRAPHY    %%%%%%%%%%%

\bibliographystyle{alpha}

\begin{thebibliography}{aaaa}



\bibitem{aus} L. Auslander, {\it Bieberbach's Theorem on Space Groups and Discrete Uniform Subgroups of Lie
Groups.} Ann. of Math. 71 (1960), 579--590.

%\bibitem{AGH}
%L. Auslander, L. Green, L.; F. Hahn. {\it Flows on homogeneous spaces} With the assistance of L. Markus and W. Massey, and an appendix by L. Greenberg. Annals of Mathematics Studies, No. 53 Princeton University Press, Princeton, N.J. 1963.


\bibitem{BG} E. Bombieri and W. Gubler, 
{\it Heights in Diophantine geometry.}
New Mathematical Monographs, 4. Cambridge University Press, Cambridge, 2006.


\bibitem{CB0}J.-P. Conze and S. Le Borgne,
{\it M\'ethode de martingales et flot g\'eod\'sique sur une surface de courbure constante n\'egative.}
Ergodic Theory Dynam. Systems 21 (2001), no. 2, 421--441.

\bibitem{CB}  J.-P. Conze and S. Le Borgne,
{\it  Le TCL pour une classe d'automorphismes non hyperboliques de vari\'et\'es.}
Preprint, 2002.

\bibitem{CG}  L. Corwin and F. Greenleaf, 
{\it Representations of nilpotent Lie groups and their applications. Part I. Basic theory and examples. }
Cambridge Studies in Advanced Mathematics, 18. Cambridge University Press, Cambridge, 1990.

%\bibitem{Dam07} Damjanovi\'{c}, D., {\it  Central extensions of simple Lie groups and rigidity of some abelian partially hyperbolic abelian actions},  Journal of Modern Dynamics,  v 1, no. 4, (2007), 665--688.

%\bibitem{DamKat05} Damjanovi\'{c}, D.; Katok, A., {\it Periodic cycle functionals and cocycle rigidity for certain partially hyperbolic $\R^k$-actions}. Discrete Contin. Dyn. Syst. 13 (2005), no. 4, 985--1005.

%\bibitem{DamKat10} Damjanovi\'{c}, D.; Katok, A., {\it  Local rigidity of partially hyperbolic actions. II. The geometric method and restrictions of Weyl chamber flows on $SL(n,,R)/\Gamma$}, Int Math, Res. Notes, 2010, to appear.

\bibitem{dekimpe0} K. Dekimpe, {\it What is \ldots an infra-nilmanifold endomorphism?}
 Notices Amer. Math. Soc. 58 (2011), no. 5, 688--689.

\bibitem{dekimpe} K. Dekimpe, {\it What an infra-nilmanifold endomorphism really should be.} arXiv:1008.4500.

\bibitem{Dolgopyat}  D. Dolgopyat,  {\it Limit theorems for partially hyperbolic systems}. Trans. Amer. Math. Soc. 356 (2004), no. 4, 1637--1689.

\bibitem{EW} G. Everest and T. Ward, 
{\it Heights of polynomials and entropy in algebraic dynamics.}
Universitext. Springer-Verlag London, Ltd., London, 1999.

\bibitem{FKS} D. Fisher, B. Kalinin and R. J. Spatzier, 
{\it Global Rigidity of  Anosov actions of Higher Rank Anosov Actions on Tori and Nilmanifolds}. 
To appear in J. Amer. Math. Soc.; arXiv:1110.0780.


\bibitem{GorodnikSpatzierII}
A. Gorodnik and R. Spatzier. {\it   Mixing Properties of Commuting Nilmanifold Automorphisms}.
Preprint.

\bibitem{Green-Tao}  B. Green and T. Tao, 
{\it The quantitative behaviour of polynomial orbits on nilmanifolds.}, Ann. of Math. (2) 175 (2012), no. 2, 465--540. 

\bibitem{HH} P. Hall and C. C. Heyde, {\it Martingale limit theory and its application}. Academic
Press, New York, 1980.

%\bibitem{Hurder} S. Hurder
%{\it  Affine Anosov actions} Michigan Math. J. 40 (1993), no. 3, 561--575.

%\bibitem{KS}  Katok, A.; Spatzier, R. J. {\it First cohomology of Anosov actions of higher rank abelian groups and applications to rigidity} IHES Publ. Math. No. 79 (1994), 131--156. 

\bibitem{Katz-bernoulli}   Y. Katznelson, {\it Ergodic automorphisms of $T^n$ are Bernoulli shifts.} Israel J. Math. 10 (1971), 186--195.

\bibitem{Kleinbock-Margulis} D. Y. Kleinbock and G. A. Margulis,
{\it Bounded orbits of nonquasiunipotent flows on homogeneous spaces.}
Sinai's Moscow Seminar on Dynamical Systems, 141--172, Amer. Math.
Soc. Transl. Ser. 2, 171, Amer. Math. Soc., Providence, RI, 1996.


\bibitem{LB98}
S. Le Borgne, {\it Un probleme de regularite dans l'equation de cobord. }   Fascicule de probabilites (Rennes, 1998), 10 pp.,
Publ. Inst. Rech. Math. Rennes, 1998, Univ. Rennes I, Rennes, 1998. 

\bibitem{LB99}
S. Le Borgne, {\it Limit theorems for non-hyperbolic automorphisms of the torus.} Israel J. Math. 109 (1999), 61--73.

\bibitem{LB02} S. Le Borgne, {\it Principes d'invariance pour les flots diagonaux sur $\hbox{\rm SL}(d,R)/\hbox{\rm
    SL}(d,Z)$.}
 Ann. Inst. H. Poincar\'e Probab. Statist.  38 (2002), no. 4, 581--612.

%\bibitem{LeBorgne-Pene}  S. Le Borgne and F. P\`{e}ne, {\it  Vitesse dans le th\'eor\'eme limite central pour certains syst\'emes dynamiques quasi-hyperboliques}. Bull. Soc. Math. France 133 (2005), no. 3, 395--417.


%\bibitem{Lee} Lee, K.B., {\it There are only finitely many infra-nilmanifolds under each nilmanifold. } Quart. J. Math. Oxford Ser. (2) 39 (1988), no. 153, 61--66.

\bibitem{Leonov} V. P. Leonov, {\it  On the central limit theorem for ergodic endomorphisms of compact
    commutative groups. }  Dokl. Akad. Nauk SSSR 135 (1960), 258--261. 

\bibitem{l} D. A. Lind,
{\it Ergodic group automorphisms and specification.}
 Ergodic theory (Proc. Conf., Math. Forschungsinst., Oberwolfach, 1978), pp. 93--104, 
Lecture Notes in Math., 729, Springer, Berlin, 1979.

\bibitem{Lind} D. A. Lind, {\it Dynamical properties of quasihyperbolic toral automorphisms.} Ergodic Theory Dynamical Systems 2 (1982), no. 1, 49--68. 


\bibitem{Livsic}  A. N. Livsic,
 {\it Cohomology of dynamical systems.}  Izv. Akad. Nauk SSSR Ser. Mat. 36 (1972), 1296--1320.
  
\bibitem{Marcuard} J.-C. Marcuard,
{\it Produits semi-directs et application aux nilvari\'{e}t\'{e}s et aux tores en th\'eorie ergodique. }
 Bull. Soc. Math. France 103 (1975), no. 3, 267--287.


\bibitem{Marcus} B. Marcus, {\it A note on periodic points for ergodic toral automorphisms. } Monatsh. Math. 89 (1980), no. 2, 121--129. 

%\bibitem{Mozes}  Mozes, S., {\it Mixing of all orders of Lie groups actions.} Invent. Math. 107 (1992), no. 2, 235--241. 

\bibitem{Parry} W. Parry,  {\it  Ergodic properties of affine transformations and flows on nilmanifolds.}  Amer. J. Math. 91 (1969), 757--771.

\bibitem{Pene}   E. P\`{e}ne,  {\it Averaging method for differential equations perturbed by dynamical systems.} ESAIM Probab. Statist. 6 (2002), 33--88.


\bibitem{RT} J. Rauch and M. Taylor, {\em Regularity of functions smooth along foliations, and elliptic
    regularity.} J. Funct. Anal. 225 (2005), no. 1, 74--93.

%\bibitem{R}  M. S. Raghunathan, Discrete subgroups of Lie groups.
%Ergebnisse der Mathematik und ihrer Grenzgebiete, Band 68. Springer-Verlag, New York-Heidelberg, 1972. 

\bibitem{Rok} V. A. Rokhlin,  {\it On the fundamental ideas of measure theory.} American 
Mathematical Society Translations, Series 1 10 (1962), 1--54.

\bibitem{Rudolph78}  D. Rudolph,
{\it Classifying the isometric extensions of a Bernoulli shift. }
J. Analyse Math. 34 (1978), 36--60. 

\bibitem{Rudolph} D. Rudolph, {\it Ergodic behaviour of Sullivan's geometric measure on a geometrically finite hyperbolic manifold.} Ergodic Theory Dynam. Systems 2 (1982), no. 3-4, 491--512.


%\bibitem{Schmidt} Schmidt, K., {\it Mixing automorphisms of compact groups and a theorem by Kurt Mahler}, 
%Pac. J. Math. Vol. 137, No. 2, 1989, 371--385.

\bibitem{sv}
P. Samek and D. Voln\'y,
{\it Uniqueness of a martingale-coboundary decomposition of stationary processes.}
Comment. Math. Univ. Carolin. 33 (1992), no. 1, 113--119.

\bibitem{S}  A. Starkov, {\it Dynamical systems on homogeneous spaces. }
Translations of Mathematical Monographs, 190. American Mathematical Society, Providence, RI, 2000.

\bibitem{Veech} W. Veech, {\it Periodic points and invariant pseudomeasures for toral endomorphisms.} Ergodic Theory Dynam. Systems 6 (1986), no. 3, 449--473. 

\bibitem{Viana}  M. Viana, {\it Stochastic Dynamics of Deterministic Systems}, Brazilian Math. Colloquium 1997, IMPA.


%\bibitem{Walters} P. Walters, {\it An introduction to ergodic theory. }
%Graduate Texts in Mathematics, 79.  Springer-Verlag, New York-Berlin, 1982

%\bibitem{Wang10} Wang, Z.J., {\it Local Rigidity of Partially Hyperbolic Actions},  JMD, vol. 4, no. 2, 2010, 271--327.

\bibitem{Wilkinson}  A. Wilkinson, {\it  The cohomological equation for partially hyperbolic
    diffeomorphisms.} To appear in Asterisque; arXiv:0809.4862.


\end{thebibliography}

\end{document}